\documentclass[12pt]{amsart}

\usepackage{latexsym,amsmath,amsthm,amssymb,amscd,amsfonts,diagrams,graphics,ogonek}

\DeclareFontFamily{U}{rsf}{}
\DeclareFontShape{U}{rsf}{m}{n}{
  <5> <6> rsfs5 <7> <8> <9> rsfs7 <10-> rsfs10}{}
\DeclareMathAlphabet{\mathscr}{U}{rsf}{m}{n}

\DeclareMathAlphabet{\mathgth}{U}{euf}{m}{n}

\DeclareFontFamily{U}{cyr}{}
\DeclareFontShape{U}{cyr}{m}{n}{
  <5> wncyr5 <6> wncyr6 <7> wncyr7 <8> wncyr8 <9> wncyr9 <10-> wncyr10}{}
\DeclareMathAlphabet{\mathcyr}{U}{cyr}{m}{n}

\input cyracc.def

\newcommand{\gAb}{\mathgth{Ab}}
\newcommand{\gSim}{\mathgth{Sim}}
\newcommand{\cA}{{\mathscr A}}
\newcommand{\cB}{{\mathscr B}}

\newcommand{\sC}{{\mathcal C}}
\newcommand{\sD}{{\mathcal D}}

\newcommand{\cE}{{\mathscr E}}
\newcommand{\cF}{{\mathscr F}}
\newcommand{\cG}{{\mathscr G}}
\newcommand{\cH}{{\mathscr H}}
\newcommand{\cI}{{\mathscr I}}

\newcommand{\cL}{{\mathscr L}}
\newcommand{\cM}{{\mathscr M}}
\newcommand{\cN}{{\mathscr N}}
\newcommand{\cO}{{\mathscr O}}
\newcommand{\cP}{{\mathscr P}}

\newcommand{\cU}{{\mathscr U}}

\newcommand{\FMYX}{\Phi_{Y\ra X}}
\newcommand{\FMXX}{\Phi_{X\ra X}}
\newcommand{\FMYY}{\Phi_{Y\ra Y}}
\newcommand{\FMXY}{\Phi_{X\ra Y}}

\newcommand{\FMYZ}{\Phi_{Y\ra Z}}
\newcommand{\FMXZ}{\Phi_{X\ra Z}}

\newcommand{\FMMX}{\Phi_{M\ra X}}

\newcommand{\sHom}{\underline{\mathrm{Hom}}}
\newcommand{\sExt}{\underline{\mathrm{Ext}}}

\newcommand{\D}{{\mathbf D}_{\mathrm{coh}}^b}
\newcommand{\KK}{{\mathbf K}}

\newcommand{\Dn}{{\mathbf D}}

\newcommand{\chk}{{\scriptscriptstyle\vee}}
\newcommand{\R}{\mathbf{R}}
\newcommand{\Ld}{\mathbf{L}}
\newcommand{\lotimes}{\stackrel{\Ld}{\otimes}}

\newcommand{\pt}{{\mathrm{pt}}}

\DeclareMathOperator{\ExFun}{ExFun}

\newcommand{\Hom}{{\mathrm{Hom}}}

\DeclareMathOperator{\Tor}{Tor}

\DeclareMathOperator{\Proj}{Proj}

\DeclareMathOperator{\ch}{ch}

\DeclareMathOperator{\id}{id}

\DeclareMathOperator{\Id}{Id}

\DeclareMathOperator{\Cl}{Cl}

\DeclareMathOperator{\Ob}{Ob}

\DeclareMathOperator{\Ext}{Ext}

\DeclareMathOperator{\Cone}{Cone}
\DeclareMathOperator{\GHilb}{{\mathit G}-Hilb}
\newcommand{\adjoint}{\dashv}

\newcommand{\ra}{\rightarrow}
\newcommand{\lra}{\longrightarrow}

\newcommand{\scdot}{{\,\cdot\,}}

\newcommand{\C}{\mathbf{C}}

\newcommand{\Q}{\mathbf{Q}}
\newcommand{\Z}{\mathbf{Z}}

\newcommand{\gCoh}{\mathgth{Coh}}

\newcommand{\iso}{\cong}

\newcommand{\pj}{\mathbf{P}}

\newcommand{\bone}{{\mathbf 1}}
\newarrow{Equal} =====
\newcommand{\del}{\partial}
\theoremstyle{plain}
\newtheorem{theorem}{Theorem}[section]
\newtheorem{lemma}[theorem]{Lemma}
\newtheorem{corollary}[theorem]{Corollary}
\newtheorem{proposition}[theorem]{Proposition}
\newtheorem{conjecture}[theorem]{Conjecture}
\theoremstyle{definition}
\newtheorem{definition}[theorem]{Definition}

\newtheorem{definition-theorem}[theorem]{Definition-Theorem}
\newtheorem{example}[theorem]{Example}

\newtheorem{exercise}[theorem]{Exercise}
\theoremstyle{remark}
\newtheorem{remark}[theorem]{Remark}

\renewcommand{\phi}{\varphi}

\begin{document}

\author{Andrei C\u ald\u araru}

\title{Derived categories of sheaves: a skimming}

\date{}

\begin{abstract}
These lecture notes were prepared for the workshop ``Algebraic
Geometry: Presentations by Young Researchers'' in Snowbird, Utah, July
2004, and for the autumn school in \L uk\k ecin, Poland, September 2004.
In six lectures I attempted to present a non-technical introduction to
derived categories of sheaves on varieties, including some important
results like:
\begin{itemize}
\item[--] the description of the derived category of $\pj^n$;
\item[--] Serre duality;
\item[--] the Bondal-Orlov result about derived categories of varieties with ample or antiample canonical class;
\item[--] the Mukai-Bondal-Orlov-Bridgeland criterion for equivalence;
\item[--] equivalences for K3 surfaces and abelian varieties;
\item[--] invariance of Hochschild homology and cohomology under derived
equivalence.
\end{itemize}
\end{abstract}

\maketitle

\tableofcontents

\section{Lecture 1: The derived category}

\subsection{}
The derived category is a rather complicated object, having its roots
in efforts by Grothendieck and Verdier to put firmer foundations on
homological algebra than those originally laid out by Grothendieck in
his T\^{o}hoku paper~\cite{Toh}.  While the need for the derived
category originally appeared as a technical tool needed for
generalizing Poincar\'e and Serre duality to relative settings, lately
(especially since Kontsevich's statement of Homological Mirror
Symmetry~\cite{Kon}) derived categories have taken a life of their
own.

\subsection{}
The motto of derived categories is~(\cite{Tho})
\begin{quote}
{\em Complexes {\bf good}, homology of complexes {\bf bad}.}
\end{quote}

\noindent
To justify this motto, we will take some cues from algebraic
topology.  (Parts of this exposition have been inspired by Richard
Thomas' article~\cite{Tho}, which is highly suggested reading.)
\vspace{2mm}

\noindent
The motivating problem is the following:
\vspace{2mm}

\noindent {\bf Problem.} Define an invariant of simplicial complexes
which allows us to decide when the topological realizations $|X|$ and
$|Y|$ of simplicial complexes $X$ and $Y$ are homotopy equivalent.
\vspace{2mm}

The first choice that comes to mind is to use homology, as
we know that if $|X|$ is homotopy equivalent to $|Y|$, then $H_i(X) =
H_i(|X|) \iso H_i(|Y|) = H_i(Y)$.  However, this turns out to be too
limited:
\vspace{2mm}

\noindent {\bf Fact.}  There exist topological spaces $X$ and $Y$
such that $H_i(X)$ is isomorphic to $H_i(Y)$ for every $i$, but $X$ is
not homotopy equivalent to $Y$.
\vspace{2mm}

This should be understood as saying that knowing the homology of a
space gives only limited information about its homotopy type.
(Homology bad.)  However, recall how $H_\scdot(X)$ is defined: it is
the homology of the chain complex $C_\scdot(X)$.  The moral that
complexes are good suggests that we should take as the invariant of
$X$ the entire chain complex $C_\scdot(X)$.  Indeed, we have the
following theorem:
\vspace{2mm}

\noindent {\bf Theorem (Whitehead).}  {\em Simplicial complexes $X$
and $Y$ have homotopy equivalent geometric realizations $|X|$ and
$|Y|$ if, and only if, there exists a simplicial complex $Z$ and
simplicial maps $f:Z \ra X$ and $g:Z\ra Y$ such that the maps $f_*:
C_\scdot(Z) \ra C_\scdot(X)$ and $g_*:C_\scdot(Z) \ra C_\scdot(Y)$ are
{\em quasi-isomorphisms}, i.e., the induced maps 
\[ f_*:H_i(Z) \ra H_i(X),\quad g_*:H_i(Z) \ra H_i(Y) \]
are isomorphisms for all $i$.}
\vspace{2mm}

What we see here is that in order to obtain a complete homotopy
invariant, we need to remember not only the homology $H_\scdot(X)$ of
the chain complex $C_\scdot(X)$, but rather the complex itself.  It is
not enough to require $H_i(X)\iso H_i(Y)$ as abstract vector spaces;
we need a chain map $C_\scdot(X) \ra C_\scdot(Y)$ inducing these
isomorphisms.

\subsection{}
Returning to algebraic geometry, recall the construction of derived
functors: sheaf cohomology, Ext, Tor, etc.\ are defined by means of
either injective or locally free resolutions.  In all these instances,
the invariant that is being defined (say, $H^i(X,\cF)$ for a sheaf
$\cF$) is obtained as the homology of a {\em complex} (in the example,
the complex $\Gamma^\scdot(X, \cI^\scdot)$ of global sections of an
injective resolution $\cI^\scdot$ of $\cF$).  The motto stated in the
beginning suggests that this is the wrong procedure: instead of
remembering just the vector spaces $H^i(X, \cF)$, we should remember
the entire complex $\Gamma^\scdot(X, \cI^\scdot)$.  This is
precisely what the derived category will do.

\subsection{The homotopy category}
Having decided that the main objects of study will be complexes (of
abelian groups, vector spaces, coherent sheaves, etc.), we'd like them
to form a category.  Indeed, we know for example that maps between
topological spaces induce chain maps between their chain complexes,
Whitehead's theorem is phrased in terms of these maps, etc.  Thus,
we'll need to study not just complexes, but also maps between them.

A natural first choice would be to start with a base abelian
category $\cA$ (abelian groups, vector spaces, sheaves, etc.) and
construct a new category whose objects are complexes of objects of
$\cA$, and maps are chain maps.  More precisely, objects are of the
form
\[ A^\scdot = \cdots \stackrel{d^{i-2}}{\lra} A^{i-1}
 \stackrel{d^{i-1}}{\lra} A^i \stackrel{d^i}{\lra} A^{i+1}
\stackrel{d^{i+1}}{\lra} \cdots, \] 
where the $A^i$'s are objects of $\cA$, $d^{i+1}\circ d^i = 0$ for all $i$,
and a morphism between $A^\scdot$ and $B^\scdot$ is a chain map
$f:A^\scdot \ra B^\scdot$, i.e., a collection of maps $f^i:A^i \ra
B^i$ such that the obvious squares commute.

Observe that at this point we have switched to cohomological notation
-- upper indices, $A^i = A_{-i}$.  From now on we will keep this
notation, and all our complexes will have their differentials increase
the degree.

Since we want to do something similar to homotopy theory, we'd like to
treat homotopic maps between topological spaces as being equal; in
particular, if $f, g: Y\ra X$ are homotopic maps between simplicial
complexes, we'd like the induced maps $f^*, g^*:C^\scdot(X) \ra
C^\scdot(Y)$ to be equal.  Now recall that a homotopy between $f$ and
$g$ yields maps $h^i:C^i(X) \ra C^{i-1}(Y)$ (which may not commute
with the differentials) such that
\[ f^i - g^i = d_Y^{i-1} \circ h^i + h^{i+1} \circ d_X^i. \]
It is an easy exercise to show that if we declare two
chain maps $f^\scdot, g^\scdot:A^\scdot \ra B^\scdot$ to be
homotopic (in the abstract category we started to define) if and only
if there exists a map $h^\scdot$ as above, then homotopy is an
equivalence relation, and composing homotopic maps yields homotopic
compositions.

\begin{definition}
Let $\cA$ be an abelian category (e.g., the category of abelian groups
$\gAb$).  We define the homotopy category of $\cA$, $\KK(\cA)$, to be
the category whose objects are complexes of objects of $\cA$, and
homomorphisms between complexes are chain maps modulo the homotopy
equivalence relation.
\end{definition}

\subsection{Injective resolutions}
Another fundamental reason why we want to treat homotopic maps as
being equal is given by the following exercise:

\begin{exercise}
Let $A, B \in \Ob\cA$ (if you feel more comfortable this way, say $A,
B$ are abelian groups, or coherent sheaves), and let $f: A\ra B$ be a
morphism.  Let $0\ra A \ra I_A^\scdot$ and $0\ra B\ra I_B^\scdot$ be
injective resolutions of $A$ and $B$, respectively.  Then prove that
that the map $f:A\ra B$ can be lifted to a map $\bar{f}$ of complexes
$\bar{f}:I_A^\scdot \ra I_B^\scdot$ such that $H^0(\bar{f}) = f$ (this
is what it means that $\bar{f}$ lifts $f$).  Argue that it is not true
that $\bar{f}$ is unique, but prove that any two liftings of $f$ are
homotopic.
\end{exercise}

The above result should be understood as saying that morphisms between
injective resolutions are, up to homotopy (i.e., in $\KK(\cA)$) the same
as morphisms between the objects being resolved:
\[ \Hom_{\KK(\cA)}(I_A^\scdot, I_B^\scdot) \iso \Hom_{\cA}(A, B). \]

\subsection{The derived category}
While the above definition provides a good framework for studying the
homotopy category of spaces, we'd like to have a setting in which
Whitehead's theorem can be naturally expressed.  More precisely,
observe that the way $\KK(\gAb)$ was constructed, there is a natural
functor
\[ C^\scdot : \gSim/\mathrm{homotopy} \lra \KK(\gAb) \]
which associates to a simplicial complex $X$ its cochain complex
$C^\scdot(X)$, and to a simplicial map the associated chain map.  (The
category $\gSim/\mathrm{homotopy}$ is the category of simplicial
complexes with morphisms simplicial maps modulo homotopy.)

The ideal result would be that $C^\scdot$ maps homotopy equivalent
simplicial complexes to isomorphic cochain complexes in $\KK(\gAb)$.
However, Whitehead's theorem involves the use of a third simplicial
complex $Z$: this is necessary because a homotopy between $|X|$ and
$|Y|$ may not respect the simplicial structure of $X$, and a further
refinement of this structure may be needed.  We notice though that for
our purposes, it is enough to pretend that chain maps that induce
isomorphisms on homology are themselves isomorphisms.  Then we can
pretend that the maps $C^\scdot(Z) \ra C^\scdot(X)$, $C^\scdot(Z) \ra
C^\scdot(Y)$ are invertible (even though most likely they are not),
and we get out of the roof
\[
\begin{diagram}[height=1.5em,width=1.5em,labelstyle=\scriptstyle]
& & C^\scdot(Z) & & \\
& \ldTo & & \rdTo & \\
C^\scdot(X) & & \pile{\lDashto\\ \rDashto} & & C^\scdot(Y)
\end{diagram}
\]
``maps'' $C^\scdot(X) \ra C^\scdot(Y)$ and $C^\scdot(X) \ra C^\scdot(Y)$
which are ``inverse'' to one another.

\begin{definition}
A chain map of complexes $f:A^\scdot \ra B^\scdot$ is called a {\em
  quasi-isomorphism} if the induced maps $H^i(f):H^i(A) \ra H^i(B)$ are
  isomorphisms for all $i$.
\end{definition}

Note that a quasi-isomorphism is often not invertible: the
map of complexes
\[ 
\begin{diagram}[height=1.5em,width=1.5em,labelstyle=\scriptstyle]
\cdots & \rTo & 0 & \rTo & \Z & \rTo^{\cdot 2} & \Z & \rTo & 0 & \rTo &
\cdots \\
& & \dTo & & \dTo & & \dTo & & \dTo & & \\
\cdots & \rTo & 0 & \rTo & 0 & \rTo & \Z/2\Z & \rTo & 0 & \rTo &
\cdots
\end{diagram}
\]
is obviously a quasi-isomorphism but does not have an inverse.  Also,
there are complexes that have the same homology, but are not
quasi-isomorphic, for a (quite non-trivial) example see
\[ \C[x, y]^{\oplus 2} \stackrel{(x,y)}{\lra} \C[x, y] \mbox{ and }
\C[x,y] \stackrel{0}{\lra} \C. \]

\subsection{}
The derived category $\Dn(\cA)$ of the initial abelian category
$\cA$ is obtained by ``pretending'' that quasi-isomorphisms in
$\KK(\cA)$ are isomorphisms.  This process is called {\em localization},
by analogy with localization of rings (where we pretend that certain
ring elements are invertible).  The analogy with rings makes it easy
to see how we should construct the localization: morphisms in
$\Dn(\cA)$ between $A^\scdot$ and $B^\scdot$ will be roofs
\[
\begin{diagram}[height=1.5em,width=1.5em,labelstyle=\scriptstyle]
& & C^\scdot & & \\
& \ldTo^f & & \rdTo^g & \\
A^\scdot & &\rDashto & & B^\scdot,
\end{diagram}
\]
with $f,g$ morphisms in  $\KK(\cA)$ and $f$ a quasi-isomorphism.  This
roof  represents  $g\circ  f^{-1}$  (despite  the  obvious  fact  that
$f^{-1}$ does  not exist) -- in  much the same way  the fraction $3/4$
represents  $3\cdot 4^{-1}$  in $\Q$,  even though  $4^{-1}$  does not
exist in  $\Z$. (The fact that only  one roof suffices, and  we do not
need arbitrarily  long zig-zags is  a non-trivial exercise.)   I won't
give all  the details  of how to  construct the localization,  but the
above  picture  should make  it  quite  clear  what the  morphisms  in
$\Dn(\cA)$ are.

It is perhaps useful to emphasize again the fact that the objects of
$\Dn(\cA)$ are just good old complexes of objects in $\cA$.
What is much less obvious in the derived category is what morphisms
are -- they are roofs of homotopy classes of chain maps, and the fact
that the morphisms from $A^\scdot$ to $B^\scdot$ may involve a third
complex $C^\scdot$ (which is not given {\em a priori}) makes it hard
to grasp what morphisms are.

\subsection{Yoneda Ext's}
An unexpected consequence of localization is the fact that the derived
category makes it easy to talk about extensions.  Say 
\[ 0 \lra A \stackrel{f}{\lra} B \stackrel{g}{\lra} C \lra 0 \]
is a short exact sequence, thought of as making $B$ an extension of
$C$ by $A$.  Then we have maps of complexes
\[ 
\begin{diagram}[height=1.25em,width=2em,labelstyle=\scriptstyle]
\cdots & \rTo & 0 & \rTo & C & \rTo & 
\cdots & \quad\quad & C\\
& & \uTo & & \uTo^g & & & & \\
\cdots & \rTo & A & \rTo^{f} & B & \rTo &
\cdots & & \dTo \\
& & \dTo & & \dTo & & & & \\
\cdots & \rTo & A & \rTo & 0 & \rTo & \cdots & & A[1]
\end{diagram}
\]
where the first map is a quasi-isomorphism.  Composing the inverse of
the first map with the second we get the map at the right, $C \ra
A[1]$, where $C$, $A$ represent the complexes which have $C$, $A$,
respectively in position 0, and zeros everywhere else, and [1] denotes
the shifting of a complex by 1 to the left.

We'll see later (when we talk about derived functors) that
\[ \Hom_{\Dn(\cA)}(C, A[1]) = \Ext^1_{\cA}(C, A), \]
the right hand side being computed the usual way, with injective
resolutions.  So we see the connection of extensions (as in exact
sequences) and elements of $\Ext$-groups.

\begin{exercise}
Generalize the above construction to extensions of arbitrary length,
and construct the Yoneda product of extensions.
\end{exercise}

\section{Lecture 2: Triangles, exactness, and derived functors}

\subsection{}
One of the most important things that we lost in passing to the
homotopy category is the ability to say that a sequence of morphisms
is exact: we no longer have notions of kernel and cokernel ($\KK(\cA)$
is not an abelian category).  Verdier's initial contribution to the
development of derived category was the observation that a form of
exactness is still preserved, in the notion of {\em exact triangles}.

An exact triangle is first off a triangle, i.e., a collection of 3
objects and 3 morphisms as follows
\[ 
\begin{diagram}[height=1.5em,width=2em,labelstyle=\scriptstyle]
A & & \rTo^f & & B \\
  & \luTo_h & & \ldTo_g & \\
&& C &&
\end{diagram}
\]
where $f$ and $g$ are regular morphisms, and $h$ shifts degree by one,
$h:C \ra A[1]$.  (The derived category will be a {\em triangulated
category}, and the data of the shift functor and exact triangles
are precisely what is needed to specify such a category.)  

Not all triangles are exact (just as not every 3 groups and 2
morphisms between them form a short exact sequence).  The way we'll
say which triangles are exact is by giving some distinguished
triangles which we declare to be exact, and by declaring any other
triangle isomorphic to a distinguished triangle to also be exact.

The distinguished triangles are built via the {\em cone construction}
(for an account of the topological background on this see~\cite{Tho}):
given a chain map of complexes $f:A^\scdot \ra B^\scdot$, let $C$ be
the complex with
\[ C^i = A^{i+1} \oplus B^i \]
and differential
\[ d_C^i = \left ( \begin{array}{cc}
d_A^{i+1} & 0 \\
f^{i+1} & d_B^i
\end{array}
\right ). 
\]
It is then clear how to define the maps $g:B^\scdot \ra C^\scdot$ and
$h:C^\scdot \ra A^{\scdot+1} = A[1]$, to build a triangle.  The object
$C^\scdot$ is called the {\em cone} of the morphism $f$, and the
resulting triangle is called a {\em distinguished triangle}.

All this is done in $\KK(\cA)$.  This endows $\KK(\cA)$ with the structure
of a triangulated category (i.e., a category with a shift functor and
a collection of exact triangles satisfying a number of complicated
axioms), and one proves that localization by quasi-isomorphisms
preserves this triangulated structure.

\begin{exercise}
Check that if $0\ra A \ra B \ra C\ra 0$ is a short exact sequence,
then 
\[ 
\begin{diagram}[height=1.5em,width=2em,labelstyle=\scriptstyle]
A & & \rTo & & B \\
  & \luTo & & \ldTo & \\
&& C &&
\end{diagram}
\]
is an exact triangle in $\Dn(A)$, even though it is not necessarily so
in $\KK(\cA)$.
\end{exercise}

\begin{remark}
One of the axioms of a triangulated category states that given a
diagram
\[
\begin{diagram}[height=1.5em,width=2em,labelstyle=\scriptstyle]
A & & \rTo & & B \\
  & \luTo & & \ldTo & \\
\dTo^f && C && \dTo^g\\
&& \dDashto && \\
A' & \rLine & \VonH & \rTo & B' \\
& \luTo & & \ldTo & \\
&& C' && 
\end{diagram}
\]
where $A$, $B$, $C$ and $A'$, $B'$, $C'$ form exact triangles, and the
morphisms $f$ and $g$ are given such that the square that they are on
two sides of commutes, then there exists a map $C\ra C'$ such that all
the squares commute.  The most serious theoretical problem of derived
categories, and which leads people to believe that in the future we
should start replacing them with more complicated things (like
dg-categories or $\cA_\infty$-categories), is the fact that this
fill-in is not unique.  In most constructions of derived categories
there exists a preferred fill-in, but the derived category data does
not specify this.
\end{remark}

\subsection{Injective resolutions and derived functors}
When dealing with sheaves we are most often also very interested in
functors: $f_*$, $\Hom$, $f^*$, $\otimes$, $\Gamma(X, -)$, etc.  Given
a complex $A^\scdot \in \Dn(\cA)$, the na\"\i ve way to apply a functor $F$
to it would be to just apply $F$ to every object and every morphism in
$A^\scdot$.  But, as we know, most functors are not exact; in derived
category language, they will not take exact triangles to exact
triangles.  

The solution to this is to look for a ``best approximation'', a
functor $\R F$ (in the case of a left exact functor; or $\Ld F$, in
the case of a right exact functor) which behaves like $F$ to some
extent, and which is exact, in the sense that it takes exact triangles
to exact triangles.

I will not make precise the sense in which $\R F$ is a best
approximation to $F$.  But, for example, the following property will
hold: if $A\in \Ob(\cA)$, then 
\[ H^0(\R F(A)) = F(A). \]
(Note that even though $A$ is thought of as a complex with just one
non-zero object, $\R F(A)$ will be a complex with many non-zero
terms.)

To construct $\R F$, we first need several results on injective
resolutions. For simplicity, we will only work with bounded complexes,
i.e., complexes $A^\scdot$ such that $H^i(A^\scdot)=0 $ for $|i|\gg
0$.  The full subcategory of these complexes inside $\KK(\cA)$ will be
denoted by $\KK^b(\cA)$, and the resulting localization of
quasi-isomorphisms is denoted by $\Dn^b(\cA)$.

\begin{proposition}
In $\KK^b(\cA)$ consider a quasi-isomorphism $f:A^\scdot \ra
B^\scdot$, with $A^\scdot$ and $B^\scdot $ bounded below complexes of
{\em injective} objects of $\cA$.  Then $f$ is invertible, i.e. there
exists a chain map $g:B^\scdot \ra A^\scdot$ such that $g\circ f$ is
homotopic to $\id_{A^\scdot}$ and $f\circ g$ is homotopic to
$id_{B^\scdot}$.
\end{proposition}

\begin{exercise}
Prove it.
\end{exercise}

\begin{proposition}
Assume that $\cA$ has enough injectives, i.e., for every object $A$ of
$\cA$ we can find an injective map $0\ra A \ra I$ with $I$ injective.
Then every bounded below complex of $\KK^b(\cA)$ is quasi-isomorphic to
a complex of injectives.
\end{proposition}

\begin{exercise}
Prove it. ({\em Hint:} Take the total complex of a double complex in
which every vertical is a resolution of the corresponding object in
the initial complex.)
\end{exercise}

\begin{corollary}
\label{cor:invloc}
The composition of the following sequence of functors is an
equivalence
\[ \KK^b(\mathrm{inj}(\cA)) \ra \KK^b(\cA) \ra \Dn^b(\cA). \]
Here, the first functor is the natural inclusion of the complexes of
injectives in all complexes, and the second functor is the
localization functor $Q$, which takes a complex to itself, and a chain
map $f:A^\scdot \ra B^\scdot$ to the roof
\[
\begin{diagram}[height=1.5em,width=1.5em,labelstyle=\scriptstyle]
& & A^\scdot & & \\
& \ldTo^\id & & \rdTo^f & \\
A^\scdot & &\rDashto & & B^\scdot.
\end{diagram}
\]
\end{corollary}

\noindent
Finally, we have
\begin{proposition}
Assume that $F:\cA \ra \cB$ is a left exact functor between abelian
categories.  Then $F$ induces an exact functor
\[ \bar{F}: \KK^-(\mathrm{inj}(\cA)) \ra \KK^-(\cB) \]
obtained by directly applying $F$ to a complex of injectives in $\cA$.
\end{proposition}

\begin{exercise}
Prove it.
\end{exercise}

\begin{definition}
Assume that $\cA$ has enough injectives, and that $F:\cA \ra \cB$ is
left exact.  Fix an inverse $Q^{-1}$ to the equivalence of categories
in Corollary~\ref{cor:invloc}, and define $\R F$ to be the composite
functor
\[ \Dn^b(\cA) \stackrel{Q^{-1}}{\lra} \KK^b(\mathrm{inj}(\cA))
\stackrel{\bar{F}}{\lra} \KK^b(\cB) \stackrel{Q}{\lra} \Dn^b(\cB). \]
\end{definition}

\subsection{}
Note what the above procedure does, in very abstract terms, in the
case when we try to compute $\R F(A)$ for a single object $A$ of
$\cA$: $I^\scdot = Q^{-1}(A)$ is a complex of injectives,
quasi-isomorphic to the complex which has just $A$ in position 0 and
zero elsewhere.  This is nothing but a good old resolution of $A$.
Then it applies $F$ to this injective resolution, and declares the
result to be $\R F(A)$.  If we were to take homology, we would get
\[ H^i(\R F(A)) = H^i(F(I^\scdot)) = \R^i F(A), \]
where the second equality is the usual definition of the $i$-th right
derived functor of $F$.  What we have gained through the abstract
procedure described above is that we can now do this for arbitrary
complexes, not just for single objects, and, more importantly, the
resulting $\R F(A^\scdot)$ is a {\em complex}, not just a bunch of
disparate objects.

\begin{exercise}
Assume that $\cA$ has enough injectives.  Prove that we have
\[ \Hom_{\Dn^b(\cA)}(A, B[i]) = \Ext_{\cA}^i(A, B), \]
for $A,B\in \Ob\cA$.  In particular, if $i<0$ then $\Hom_{\Dn(\cA)}(A,
B[i]) = 0$, so morphisms in the derived category can never go ``to the
right''.
\end{exercise}

\subsection{}
What we have described above works well when $\cA$ has enough injectives
and we are dealing with a left exact functor (e.g., $\Gamma(X, -)$,
$\Hom(\cF,-)$, $f_*$).  When $F$ is right exact, we would need to have
enough projectives, and one of the main problems of sheaf theory is
that while there are enough injectives, there are almost never enough
projectives.

To get around this problem, we replace injective resolutions by
resolutions by {\em acyclic} objects.  A class of objects in $\cA$ is
called acyclic for a right exact functor $F$ if $F$ takes exact
complexes of objects in this class to exact complexes.  As an example, for
$f^*$ and $\otimes$, locally free sheaves are acyclic.  (More
generally, since there exist schemes with few locally free sheaves, we
may need to use flat sheaves, of which there are always plenty.)

If our class of acyclic objects is rich enough to be able to arrange
that every complex in $\KK(\cA)$ is quasi-isomorphic to a complex of
acyclic objects, we can repeat the construction we have used for
left-exact functors, to define $\Ld F$, the left derived functor of a
right exact functor $F$.  It will again be exact.

\subsection{}
These constructions enable us to define derived functors $\R\Gamma(X,
-)$, $\R\Hom(-, -)$, $\R\sHom(-,-)$ (these two are right-derived in
the second variable), derived tensor product, and for a morphism $f:X \ra
Y$ of schemes, $\R f_*$ and $\Ld f^*$.  The resulting complexes
compute the usual derived functors: sheaf cohomology $H^*(X, -)$,
$\Ext^*_X(-, -)$, $\sExt^*_X(-, -)$, $\Tor_*^X$, $R^* f_*$.  While the
derived pull-back will turn out to be important, it is usually not
emphasized in day-to-day algebraic geometry.

\subsection{}
The great technical advantage of using the derived category is the
fact that we can now prove that, under mild hypotheses, we have
\[ \R(F \circ G) \iso \R(F) \circ \R(G), \]
and similar relations hold for all combinations of left and right
derived functors.  Under normal circumstances, this requires the use
of a spectral sequence: for example, the Leray spectral sequence for a
morphism $f:X\ra Y$ and a sheaf $\cF$ on $X$
\[ H^i(Y, \R^j f_* \cF)  \Rightarrow H^{i+j}(X, \cF) \]
is an application of the above equality for the composition of left
exact functors
\[ \Gamma(X, -) = \Gamma(Y, f_*(-)). \]
Similarly, the local to global spectral sequence,
\[ H^i(X, \sExt^j_X(\cF, \cG)) \Rightarrow \Ext^{i+j}_X(\cF, \cG), \]
follows from
\[ \Gamma(X, \sHom_X(-,-)) = \Hom_X(-,-). \]

The reason we are able to do this is the fact that our derived
functors yield {\em complexes}; just knowing their homology would not
be enough to construct the spectral sequence, for we would not know
how to build the maps in the spectral sequence.  But having the entire
complex turns out to suffice.

Of course, if we were to explicitly compute, we would still have to
use the spectral sequences.  But using the language of derived
categories makes for much cleaner and easier statements, which would
otherwise involve large numbers of spectral sequences.

\subsection{}
\label{subsec:projdual}
Another great advantage of using the derived category is the fact that
many results that only hold for locally free sheaves in the usual
setting, now hold in general.  For example, if we denote by $\D(X)$
the bounded derived category of coherent sheaves on $X$, we have:
\begin{itemize}
\item[--] the projection formula: if $f:X\ra Y$ is a morphism of
varieties, and if $\cE\in \D(X)$, $\cF\in \D(Y)$, then we have
\[ \R f_*(\cE \lotimes \Ld f^* \cF) \iso \R f_*\cE \lotimes \cF. \]
\item[--] dualization: for $\cE, \cF, \cG\in \D(X)$ we have
\[ \R\Hom_X(\cE, \cF\lotimes \cG) \iso \R\Hom_X(\cE\lotimes \cF^\chk,
\cG), \]
where 
\[ \cF^\chk = \R\sHom_X(\cF, \cO_X). \]
\item[--] flat base change: if we are given a cartesian diagram
\[ \begin{diagram}[height=2em,width=2em,labelstyle=\scriptstyle]
X\times_Y Z & \rTo^v & X \\
\dTo_g & & \dTo_f \\
Z & \rTo^u & Y
\end{diagram} \]
with $u$ flat, then there is a natural isomorphism of functors 
\[ u^* \circ \R f_* \iso \R g_* \circ v^*. \]
\end{itemize}

\subsection{}
In all our discussions until now we have skated over several (hard)
technical details.  One of the problems has to do with bounded or
unbounded complexes.  Most of the results we stated hold true for
bounded complexes, but we need unbounded complexes (at least in one
direction) to talk about injective resolutions.  Usually this is dealt
with by working with complexes that may be infinite in one direction,
and whose homology is bounded; also, since we want to deal with
coherent sheaves, we require the homology sheaves of all complexes to
be coherent (even though we may need to use non-coherent sheaves in
some resolutions, say flat ones).  This leads to the definition of the
bounded derived category of coherent sheaves, $\D$.  From now on we'll
work with this more restrictive version of the derived category of a
scheme.

Second, if the scheme we are dealing with is not smooth, (even
though it may be projective), it is not true that every coherent sheaf
admits a {\em finite} locally free resolution.  This will cause
problems with some results, like the dualization result
in~(\ref{subsec:projdual}).  Therefore, for most of the following
results, we'll work with smooth projective schemes, unless specified
otherwise.

\section{Lecture 3: The derived category of $\pj^n$ and orthogonal 
decompositions}

\subsection{}
In today's lecture we'll cover two topics: first, we'll compute as
explicitly as possible the simplest derived category of a scheme:
$\D(\pj^n)$.  We will also skim the subject of semiorthogonal
decompositions, and talk briefly about flops and derived categories.

\subsection{}
We will give a description of $\D(\pj^n)$ in terms of generators and
relations.  The result is

\begin{theorem}
\label{thm:pn}
The derived category $\D(\pj^n)$ is generated by the exceptional
sequence
\[ \D(\pj^n) = \langle \cO(-n), \cO(-n+1), \ldots, \cO(-1), \cO
\rangle. \]
\end{theorem}

To understand this result, we first need to make clear what is meant
by ``generated.''  The idea is that the derived category has two
fundamental operations built in: shifting (a complex to the left or
right by an arbitrary amount) and cones (given objects $A^\scdot$ and
$B^\scdot$, we are allowed to pick an arbitrary morphism $f:A^\scdot
\ra B^\scdot$ and construct the cone $\Cone(f)$).  To say that
$\D(\pj^n)$ is generated by some set of objects $S$ is to say that the
smallest full subcategory of $\D(\pj^n)$ that contains the objects in
$S$ and is closed under shifts and taking cones is a subcategory
equivalent to all of $\D(\pj^n)$.

In other words, we are allowed to take objects from $S$, shift them,
take arbitrary morphisms between them, take cones of these morphisms,
and again repeat these procedures.  Up to isomorphism, every object of
$\D(\pj^n)$ should be reached after applying these operations a finite
number of times.

\subsection{}
The fact that 
\[ S = \{ \cO(-n), \cO(-n+1), \ldots, \cO(-1), \cO \} \]
generate $\D(\pj^n)$ follows from the following result, known as
Be\u\i linson's resolution of the diagonal~(\cite{Bei}).  To state it,
introduce the following notation: if $X$ and $Y$ are varieties, and
$\cE \in \D(X)$, $\cF \in \D(Y)$, then define
\[ \cE \boxtimes \cF = \pi_X^*\cE \lotimes \pi_Y^* \cF \in \D(X\times
Y). \]

\begin{proposition}[Be\u\i linson]
The following is a locally free resolution of $\cO_\Delta$ on
$\pj^n\times \pj^n$:
\begin{align*}
0 \ra \cO(-n) \boxtimes \Omega^n(n) & \ra \cO(-n+1) \boxtimes
\Omega^{n-1}(n-1) \ra \cdots \\
& \cdots \ra \cO(-1) \boxtimes \Omega^1(1) \ra \cO
\boxtimes \cO \ra \cO_\Delta\ra 0. 
\end{align*}
\end{proposition}

\begin{proof}
Fix a basis $y_0, \ldots, y_n$ of $H^0(\pj^n, \cO(1))$.  Consider the
(twist by -1 of the) Euler exact sequence of vector bundles on $\pj^n$,
\[ 0 \ra \cO(-1) \ra \cO^{n+1} \ra T(-1) \ra 0, \]
where $T$ denotes the tangent bundle of $\pj^n$.  Taking global
sections we get an isomorphism $H^0(\pj^n, \cO^{n+1}) \iso H^0(\pj^n,
T(-1))$.  A basis of $H^0(\pj^n, \cO^{n+1})$ is given by a dual basis
$y_0^\chk, \ldots y_n^\chk$ to $y_0,\ldots, y_n$, and denote by
$\del/\del y_i$ the image of $y_i^\chk$ in $H^0(\pj^n, T(-1))$.

Consider the global section $s$ of $\cO(1) \boxtimes T(-1)$ on $\pj^n
\times \pj^n$ given by
\[ s = \sum_{i=0}^n x_i \boxtimes \frac{\del}{\del y_i}, \]
where the $x_i$'s and $y_i$'s are coordinates on the first and second
$\pj^n$, respectively.  The claim is that the zeros of $s$ are
precisely along the diagonal of $\pj^n\times \pj^n$.

We check this in one coordinate patch of $\pj^n\times \pj^n$, namely
when $x_0 \neq 0, y_0\neq 0$.  (The check for the other patches is
entirely similar.)  In the patch where $y_0 \neq 0$ define affine
coordinates by $Y_i = y_i/y_0$ for $1\leq i \leq n$.  Then $\del/\del
Y_i$ for $1\leq i \leq n$ is a basis for $T$ at each point in this
patch.  From $y_i = Y_i y_0$ it follows that
\[ dY_i = \frac{y_0 dy_i + y_i dy_0}{y_0^2}. \]
Writing 
\[ \frac{\del}{\del y_i} = \sum_{i=1}^n f_i \frac{\del}{\del
  Y_i}\mbox{ with }  f_i = dY_i(\frac{\del}{\del y_i}), \]
it follows that 
\[ \frac{\del}{\del y_i} = \frac{1}{y_0}\frac{\del}{\del Y_i} \]
if $i\neq 0$, and
\[ \frac{\del}{\del y_0} = - \sum_{i=1}^n \frac{y_i}{y_0^2}
\frac{\del}{\del Y_i}. \]

Then we have
\begin{align*}
s & = \sum_{i=0}^n x_i \boxtimes \frac{\del}{\del y_i} = \sum_{i=1}^n
\frac{x_i}{y_0} \frac{\del}{\del Y_i} - \sum_{i=1}^n \frac{x_0
  y_i}{y_0^2} \frac{\del}{\del Y_i} \\
 & = \sum_{i=1}^n \frac{x_i y_0 - x_0 y_i}{y_0^2} \frac{\del}{\del
  Y_i}.
\end{align*}
Thus $s=0$ precisely when 
\[ \frac{x_i}{x_0} = \frac{y_i}{y_0} \]
for all $i$, i.e., along the diagonal of $\pj^n\times \pj^n$.

Taking the Koszul resolution for the section $s$ we get the
result.
\end{proof}

\subsection{}
\label{subsec:beiltri}
Split off the Be\u\i linson resolution into short exact sequences
\[
\begin{diagram}[width=2em,height=1.5em]
0 & \rTo & \cO(-n) \boxtimes \Omega^n(n) &\rTo& \cO(-n+1) \boxtimes \Omega^{n-1}(n-1)
&\rTo& C_{n-1} &\rTo& 0 \\ 
0 &\rTo& C_{n-1} &\rTo& \cO(-n+2) \boxtimes \Omega^{n-2}(n-2)
&\rTo& C_{n-2} &\rTo& 0 \\ 
& &&&\vdots&&&& \\
0  &\rTo& C_1 &\rTo& \cO\boxtimes \cO &\rTo& \cO_\Delta &\rTo& 0.
\end{diagram}
\]
Since short exact sequences can be viewed as particular examples of
triangles, we conclude that $\cO_\Delta$ can be generated  from
\[ \{\cO(-n) \boxtimes \Omega^n(n), \cO(-n+1) \boxtimes \Omega^{n-1}(n-1),
\ldots, \cO\boxtimes \cO \} \] 
by using $n$ triangles (on $\pj^n\times \pj^n$).

\subsection{}
The following concept will be extremely useful in the future, so we
give it a name:
\begin{definition}
\label{def:intfun}
If $X$ and $Y$ are varieties, and $\cE$ is an object in $\D(X\times
Y)$, define the {\em integral transform} with kernel $\cE$ to be
the functor
\[ \FMXY^\cE: \D(X) \ra \D(Y), \quad\quad \FMXY^\cE(-) =
\R\pi_{Y,*}(\Ld \pi_X^*(-) \lotimes \cE). \]
\end{definition}

Since all the functors involved are exact, it is obvious that if we
have an exact triangle
\[ \cE \ra \cF \ra \cG \ra \cE[1] \]
of objects on $X\times Y$, then for every $\cA\in \D(X)$ we get an
exact triangle in $\D(Y)$
\[ \FMXY^\cE(\cA) \ra \FMXY^\cF(\cA) \ra \FMXY^\cG(\cA) \ra
\FMXY^\cE(\cA)[1]. \] 

\subsection{}
Let $\Phi^\cE$ denote the integral transform with kernel
$\cE\in\D(\pj^n \times \pj^n)$, thought of as going from the {\em
  second} $\pj^n$ to the first.

Using the triangles on $\pj^n\times \pj^n$ constructed
in~(\ref{subsec:beiltri}), we conclude that for every $\cA\in
\D(\pj^n)$ there is a sequence of exact triangles on $\pj^n$
\[ 
\begin{diagram}[width=2em]
\Phi^{\cO(-n)\boxtimes\Omega^n(n)}(\cA) &\rTo
& \Phi^{\cO(-n+1)\boxtimes\Omega^{n-1}(n-1)}(\cA) &\rTo& 
\Phi^{C_{n-1}}(\cA) &\rTo& 
\Phi^{\cO(-n)\boxtimes\Omega^n(n)}(\cA)[1] \\
\Phi^{C_{n-1}}(\cA) &\rTo &
\Phi^{\cO(-n+1)\boxtimes\Omega^{n-1}(n-1)}(\cA) &\rTo& 
\Phi^{C_{n-2}}(\cA) &\rTo& 
\Phi^{C_{n-1}}(\cA)[1] \\
& &&\vdots&&& \\
\Phi^{C_1}(\cA) &\rTo& \Phi^{\cO\boxtimes \cO}(\cA) &\rTo&
\Phi^{\cO_\Delta}(\cA) &\rTo& \Phi^{C_1}(\cA)[1].
\end{diagram}
\]
Therefore we conclude that $\Phi^{\cO_\Delta}(\cA)$
is generated by
\[ \{ \Phi^{\cO(-n)\boxtimes\Omega^n(n)}(\cA),
\Phi^{\cO(-n+1)\boxtimes\Omega^{n-1}(n-1)}(\cA), \ldots,
\Phi^{\cO\boxtimes \cO}(\cA) \}. \]

\subsection{} 
Using the projection formula, it is an easy exercise to see that
\[ \Phi^{\cO_\Delta}(\cA) \iso \cA, \]
so to conclude that 
\[ \{ \cO(-n), \cO(-n+1),\ldots, \cO \} \]
generate $\D(\pj^n)$ it will suffice to argue that 
\[ \Phi^{\cO(-i)\boxtimes\Omega^i(i)}(\cA) \]
is in the subcategory generated by $\cO(-i)$.  Using the projection
formula again, we see that 
\[ \Phi^{\cO(-i)\boxtimes\Omega^i(i)}(\cA) \iso \cO(-i) \otimes_\C
\R\Gamma(\pj^n, \cA\otimes \Omega^i(i)). \]
In other words, $\Phi^{\cO(-i)\boxtimes\Omega^i}(\cA)$ is
isomorphic to a complex which has all differentials zero, and in
position $k$ it has 
\[ \dim \R^k\Gamma(\pj^n, \cA \otimes \Omega^i) \]
copies of $\cO(-i)$.  This complex is obviously generated by $\cO(-i)$
(it is a finite direct sum, i.e., cones on zero morphisms, of shifts
of $\cO(-i)$).

We conclude that every object in $\D(\pj^n)$ can be obtained by taking
at most $n$ cones on objects that look like finite direct sums of
$\cO(-i)[j]$, $-n\leq i\leq 0$, $j\in\Z$.  (And in fact, the above
proof is constructive, giving an algorithm for building the cones.
The actual algorithm is usually packaged as a spectral sequence, the
{\em Be\u\i linson spectral sequence}).

\subsection{}
It is also worth observing that 
\[ \{ \Omega^n(n), \Omega^{n-1}(n-1), \ldots, \Omega^1(1), \cO \} \]
also form a generating set for $\D(\pj^n)$.  The exact same argument
holds, but we apply the integral transforms from the first factor to
the second one.

\subsection{}
Now we turn to the second statement of Theorem~\ref{thm:pn}, namely
that 
\[ \langle \cO(-n), \cO(-n+1), \ldots, \cO \rangle \]
forms an exceptional sequence.  This is easily checked once we know
the following definition:
\begin{definition}
A sequence of objects 
\[ \langle A_n, A_{n-1}, \ldots, A_0 \rangle \]
in a triangulated category is called {\em exceptional} if
\begin{align*}
&\Ext^i(A_p, A_q)  = 0 \mbox{ for all $i$, if } p < q, \\
\intertext{and}
& \Ext^i(A_p, A_p) = \left \{ \begin{array}{ll} 0 & \mbox{ if }
  i>0,\\ \C & \mbox{ if } i = 0. \end{array} \right .
\end{align*}
(Here we wrote $\C$, if we were working over a different field than the
complex numbers this would be the ground field.)
\end{definition}

\subsection{}
The reason exceptional sequences are useful is because a derived
category constructed using one can be thought of as being pieced
together in the simplest possible way: the category generated by every
object in the exceptional sequence is equivalent to the derived
category of vector spaces, and we build the larger category by adding
these ``one dimensional subspaces'' one after another.  Furthermore,
in the new category, the smaller ones,
\[ C_i = \langle A_i, A_{i-1}, \ldots, A_0 \rangle \]
are semi-direct summands in the sense that 
\[ \Hom(C_i, \langle A_j\rangle) = 0 \mbox{ for } i<j. \]
(I.e., Hom from any object in $C_i$ to any object in the category
generated by $A_j$ is zero.)

\subsection{}
We finish with some remarks on how the calculation of $\D(\pj^n)$ can
be generalized.  Orlov, in~\cite{Orlbun}, argued that the relative
situation can be described by saying that if $\cE$ is a vector bundle
of rank $n+1$ on a projective variety $X$, then the derived category
of $\pj = \pj(\cE)$ has a semiorthogonal decomposition
\[ \D(\pj) \iso \langle \cO_\pj(-n), \cO_\pj(-n+1), \ldots,
\cO_\pj(-1), \D(X)\rangle, \]
where $\D(X)$ is seen as a full subcategory of $\D(\pj)$ via the
pull-back $p^*$ of the canonical projection $p:\pj\ra X$.

He and Bondal~\cite{BonOrl} also argued that if $Y$ is a smooth
subvariety of codimension $n+1$ of the smooth variety $X$, and if
$Z\ra X$ is the blow-up of $X$ along $Y$ with exceptional divisor $E$
(which is a projective bundle over of rank $n$ over $Y$), then
\[ \D(Z) \iso \langle \D(Y)_{-n}, \D(Y)_{-n+1}, \ldots, \D(Y)_{-1},
\D(X)\rangle. \]
Here, $\D(X)$ is viewed as a full subcategory of $\D(Z)$ via
pull-back, and $\D(Y)_{-k}$ is the full subcategory of $\D(Z)$ of
objects of the form $i_*(p^* \cE \otimes \cO_E(-k))$, where $\cE\in \D(Y)$,
$p:E\ra Y$ is the natural projection, and $i:E\ra Z$ is the inclusion.

\subsection{}
An interesting application of this was worked out by Bondal and Orlov
in~\cite{BonOrl}: say $X^-$ is a smooth threefold, and $C^-$ is a
smooth $\pj^1$ inside $X$ with normal bundle $\cO(-1)\oplus \cO(-1)$.
Blowing $C^-$ up we get a total space $Y$ and exceptional divisor $E$
isomorphic to $C^- \times \pj^1 = C^-\times C^+$.  The $C^-$ direction
in $E$ can be contracted down to get another threefold $X^+$ with
another curve $C^+$ that has the same properties inside $X^+$ as $C^-$
has inside $X^-$.  (This is the simplest example of a {\em flop}.)
While $X^-$ and $X^+$ are birational, they are in general not
isomorphic.

Bondal and Orlov showed that from the semiorthogonal decompositions
above one can argue that
\[ \D(X^-) \iso \D(X^+) \]
and they were led to conjecture that this would be true in the case
of an arbitrary flop.  Their conjecture was proven by
Bridgeland~\cite{BriFlops} and later generalized by Chen~\cite{Chen}
and Kawamata~\cite{Kaw}.  There are some indications that ideas from
derived category theory could provide great simplifications to the
minimal model program in higher dimensions, although much still needs
to be done.

\section{Lecture 4: Serre duality and ampleness of the canonical class}

\subsection{}
One of the fundamental facts that was known already to the Italian
school is that, without additional information, on a smooth
variety $X$ there are two distinguished line bundles: the trivial
bundle $\cO_X$, and the {\em canonical bundle} $\omega_X$.  It was
Serre who made more precise the r\^ole of $\omega_X$, through Serre
duality: for a coherent sheaf $\cF$ on a smooth projective variety we
have
\[ H^i(X, \cF)^\chk \iso \Ext_X^{n-i}(\cF, \omega_X) \]
for all $i$.

\subsection{}
It turns out that the canonical line bundle and Serre duality play a
crucial part in the study of derived categories.  In fact, the
historical relationship is quite the opposite: derived categories were
introduced (by Verdier and Grothendieck) in order to generalize Serre
duality to a relative context.  In today's lecture we shall review
Grothendieck's formulation of relative Serre duality, and show how it
follows easily in the case of projective morphisms between smooth
varieties from the idea of Bondal and Kapranov of a {\em Serre
functor}.  

The second half of today's lecture will cover the following question:
how much of the variety can be recovered from its derived category?
In the extreme cases, namely when $\omega_X$ is either ample or
antiample, we'll see that the variety can be completely recovered from
$\D(X)$.  Finally, even when $\omega_X$ is somewhat trivial, we argue
that the pluricanonical ring of $X$,
\[ R(X) = \bigoplus_i H^0(X, \omega_X^{\otimes i}) \]
is an invariant of $\D(X)$.

\subsection{}
Grothendieck's approach to Serre duality was motivated by the search
for a right adjoint to the push-forward functor $f_*$ for a morphism
$f:X\ra Y$ of schemes.

Recall that two functors 
\[ F:\sC \ra \sD,\quad  G:\sD \ra \sC \]
are said to be adjoint (written as $F\adjoint G$, and then $F$ is a
{\em left} adjoint to $G$, and $G$ is a {\em right} adjoint to $F$) if
there are isomorphisms
\[ \Hom_\sD(Fx, y) \iso \Hom_\sC(x, Gy), \]
natural in both variables, for every $x\in \Ob\sC$, $y\in \Ob\sD$.
For example, for a projective morphism $f:X\ra Y$, $f^*:\gCoh(Y) \ra
\gCoh(X)$, $f_*:\gCoh(X) \ra \gCoh(Y)$ are adjoint, $f^* \adjoint
f_*$.

At first, the search for a right adjoint to $f_*:\gCoh(X) \ra
\gCoh(Y)$ seems to be doomed for failure, as the following proposition
proves:
\begin{proposition}
Assume that $F\adjoint G$, and that $\sC$, $\sD$ are abelian
categories.  Then $F$ must be right exact, and $G$ must be left exact.
\end{proposition}

\begin{exercise}
Prove it!
\end{exercise}

This shows immediately that $f_*$, being left exact, can have a left
adjoint ($f^*$), but not a right adjoint (unless $f$ is something like
a closed immersion, when $f_*$ is exact).

\subsection{}
The problem, however, appears to be with the fact that we are dealing
with abelian categories, which are a bit too restrictive.  The point
is that in passing to the derived categories, we have replaced a left
exact functor $f_*$ by its right derived version $\R f_*$, which is
now exact.  So in principle there may be hope to find a right adjoint
to $\R f_*:\D(X) \ra \D(Y)$.

Let $f:X\ra \pt$ be the structure morphism of a smooth projective
variety $X$.  Let's look at what finding a right adjoint to $f_*$
would mean.  Finding a right adjoint $f^!$ to $\R f_*$ means we must have
\[ \Hom_\pt(\R f_* \cF, \cG) \iso \Hom_X(\cF, f^! \cG). \]
for $\cF\in \D(X)$, $\cG \in \D(\pt)$.  Now objects in $\D(\pt)$ are
not very interesting, they all decompose into direct sums of shifts of
the one-dimensional vector space $\C$.  (Again, we use $\C$ for our
base field, but our discussion works over an arbitrary ground field.)
So take $\cG = \C$.  Also, for simplicity, let $\cF$ consist of just a
single coherent sheaf, shifted left by $i$, which we denote by
$\cF[i]$.

Now it is easy to see that 
\[ \Hom_\pt(\R f_* \cF[i], \C) \iso H^i(X, \cF)^\chk. \]
Serre duality predicts
\[ H^i(X, \cF)^\chk \iso \Ext_X^{n-i}(\cF, \omega_X) =
\Hom_{\D(X)}(\cF[i], \omega_X[n]), \]
where $\omega_X[n]$ denotes the complex having $\omega_X$ in position
$-n$ and zero elsewhere.

Thus, if we set $f^!\C = \omega_X[n]$, Serre duality yields
\[ \Hom_\pt(\R f_* \cF[i], \C) \iso \Hom_X(\cF[i], f^! \C), \]
which is precisely the adjunction formula we want.  Making $f^!$
commute with direct sums and shifts we get a well defined $f^!:\D(\pt)
\ra \D(X)$ which is a right adjoint (in the derived sense) to $\R f_*$.

\subsection{}
We will generalize the above calculation to make it work for a much
more general class of maps $f$.  Our final result will be:
\begin{theorem}
\label{thm:gd}
Let $f:X\ra Y$ be a morphism of smooth projective schemes.  Then a
right adjoint to $\R f_*:\D(X) \ra \D(Y)$ exists, and is given by
$f^!:\D(Y) \ra \D(X)$,
\[ f^! (-) = \Ld f^*(- \otimes \omega_Y^{-1}[-\dim Y]) \otimes
\omega_X[\dim X]. \]
\end{theorem}

Note that the previous calculation that $f^! \C = \omega_X[\dim X]$
for the structure morphism $f:X\ra\pt$ of a smooth scheme $X$ agrees
with this theorem.

It is worth mentioning that much more general versions of this
statement have been proven.  Essentially, all that is required is that
$f$ be proper.  See~\cite{HarRD} for more on this.

\subsection{}
Before we begin our proof of Theorem~\ref{thm:gd} we divert for a
minute to state a different formulation of Serre duality, due to
Bondal-Kapranov~\cite{BonKap}.  

Let $X$ be a smooth, projective scheme of dimension $n$, and let $\cE,
\cF$ be vector bundles on $X$.  Then we have
\[ \Ext^i_X(\cE, \cF) \iso H^i(X, \cE^\chk\otimes \cF) \iso H^{n-i}(X,
\cE \otimes \cF^\chk \otimes \omega_X)^\chk \iso \Ext^{n-i}_X(\cF, \cE
\otimes \omega_X)^\chk. \]
Using shifts in the derived category, we can rewrite this as
\[ \Hom_X(\cE, \cF[i]) \iso \Hom_X(\cF[i], \cE \otimes
\omega_X[n])^\chk. \]
In fact, in the derived category the requirement that $\cE$ and $\cF$
be vector bundles is superfluous, so defining the functor 
\[ S_X:\D(X) \ra \D(X), \quad S_X\cE = \cE \otimes \omega_X[\dim X]
\]
we have
\[ \Hom_X(\cE, \cF) \iso \Hom_X(\cF, S_X\cE)^\chk \]
for any $\cE, \cF\in \D(X)$.

This leads Bondal and Kapranov to make the following definition:
\begin{definition}
Let $\sC$ be a $k$-linear category (i.e., $\Hom$ spaces are vector
spaces over the field $k$).  Then an equivalence $S:\sC\ra \sC$ is
called a Serre functor for $\sC$ if there exist natural, bifunctorial
isomorphisms 
\[ \phi_{A,B}:\Hom_\sC(A,B) \stackrel{\sim}{\lra} \Hom_\sC(B,
SA)^\chk \] 
for every $A,B\in\Ob(\sC)$.  (The actual isomorphisms $\phi_{A,B}$ are
part of the data required to specify the Serre functor $S$.)
\end{definition}

\subsection{}
The crucial properties of a Serre functor are summarized in the
theorem below:
\begin{theorem}[\cite{BonKap}]
Let $\sC$ be a triangulated, $k$-linear category.  Then a Serre
functor, if it exists, is unique up to a unique isomorphism of
functors.  It is exact (i.e., takes exact triangles to exact
triangles), commutes with translation, and if $F:\sC\ra \sD$ is a
triangulated equivalence, then there is a canonical isomorphism
\[ F\circ S_\sC \iso S_\sD \circ F. \]
\end{theorem}

To have a flavor of how these statements are proven, we prove that
a Serre functor must commute with translations:
\begin{align*}
\Hom(A, S(B[1])) & \iso \Hom(B[1], A)^\chk \iso \Hom(B, A[-1])^\chk \\
& \iso \Hom(A[-1], SB) \iso \Hom(A, (SB)[1]),
\end{align*}
and the isomorphism of functors 
\[ S \circ [1] \iso [1] \circ S \]
follows now from a representability theorem.

\subsection{}
It is crucial to emphasize at this point that it is {\bf not} true
that an equivalence of triangulated categories $F:\D(X)\ra \D(Y)$
between the derived categories of smooth varieties will satisfy
\[ F(\omega_X[\dim X]) \iso \omega_Y[\dim Y]. \]
The reason for this apparent mismatch is the fact that in general
derived equivalences need not take tensor products to tensor products;
thus the peculiar fact that $S_X$ is given by tensoring with
$\omega_X[\dim X]$ will not translate to $F$ mapping $\omega_X$ to
$\omega_Y$.  What is true is that $F$ will commute with the
corresponding Serre functors.  An explicit example of an equivalence
$F$ such that $F(\omega_X) \not\iso \omega_Y$ is given in Lecture 5.

\subsection{}
The great thing about Serre functors is that it allows one to convert
from a left adjoint to a right adjoint and vice versa.  Specifically,
we have
\begin{theorem}
\label{thm:leftrightadj}
Let $F:\sC\ra \sD$ be a functor between the $k$-linear categories
$\sC$, $\sD$ that admit Serre functors $S_\sC$, $S_\sD$.  Assume that
$F$ has a left adjoint $G:\sD\ra\sC$, $G\adjoint F$.  Then 
\[ H = S_\sC \circ G \circ S_{\sD}^{-1}:\sD\ra \sC \]
is a right adjoint to $F$. 
\end{theorem}

\begin{proof}
\begin{align*}
\Hom_\sD(Fx, y) & \iso \Hom_\sD(S_{\sD}^{-1}y, Fx)^\chk \iso
\Hom_\sC(GS_{\sD}^{-1}y, x)^\chk \\ & \iso \Hom_\sC(x,
S_{\sC}GS_{\sD}^{-1}y) = \Hom_\sC(x, Hy).
\end{align*}
\end{proof}

This immediately yields Theorem~\ref{thm:gd}.  One final point to note
about $f^!$ is that 
\[ \Ld  f^*(-\otimes \omega_Y[\dim Y]) \iso \Ld f^*
(-) \otimes f^*(\omega_Y)[\dim Y]. \]
Therefore we can rewrite 
\[ f^!(-) = \Ld f^*(-) \otimes \omega_{X/Y}[\dim f], \]
where 
\[ \omega_{X/Y} = \omega_X \otimes f^* \omega_Y^{-1} \]
is the relative dualizing sheaf, and 
\[ \dim f = \dim X - \dim Y. \]

\subsection{}
We now turn to the following fundamental question: how much of the
space $X$ can we recover by knowing just $\D(X)$?  The idea behind
this question comes from the following basic facts:
\begin{itemize}
\item[--] if $X$ and $Y$ are projective varieties, then $\gCoh(X) \iso
  \gCoh(Y)$ implies that $X$ is isomorphic to $Y$; therefore the
  abelian category of coherent sheaves on a scheme determines the
  scheme itself.

\item[--] similarly, if $R$ and $S$ are commutative rings that are
  Morita equivalent (their categories of modules are equivalent), then
  $R$ is isomorphic to $S$.

\item[--] there are examples of smooth, projective
  varieties $X$ and $Y$ such that $\D(X)\iso \D(Y)$ as triangulated
  categories, but $X \not \iso Y$.
\end{itemize}

Therefore it appears that the derived category is in a sense a
significantly weaker invariant than $\gCoh(X)$.  In fact, this is
precisely what makes derived categories so appealing: they are
sufficiently lax to allow for interesting pairs of varieties to be
derived equivalent, but not so lax that most everything will be
equivalent.

\subsection{}
The philosophical point here turns out to be that the farther away the
canonical class is from being trivial, the more of the variety we can
reconstruct.  More precisely, we have
\begin{theorem}[{Reconstruction -- Bondal-Orlov~\cite{BonOrl}}]
Let $X$ be a projective variety such that $\omega_X$ is either ample
or anti-ample, and let $Y$ be any projective variety.  If $\D(X) \iso
\D(Y)$ as triangulated categories, then $X\iso Y$.
\end{theorem}

In fact, the theorem can be made even more precise: it also lists what
kind of autoequivalences $\D(X)$ can have if $X$ has ample or
anti-ample canonical bundle.  Namely, only ``boring'' automorphisms
can occur -- translations, tensoring by line bundles, or pull-backs by
automorphisms of $X$ itself.

\subsection{}
We give a complete account of the reconstruction theorem
from~\cite{BonOrl}, modulo two easy results which can be found in
Bondal-Orlov~[loc.cit.].  

We begin with the observation that one can give a purely categorical
description of what elements of $\D(X)$ are structure sheaves of
points, assuming $X$ has ample canonical or anticanonical bundle:
\begin{proposition}[\cite{BonOrl}]
\label{prop:po}
Let $X$ be a smooth, projective variety with ample canonical or
anticanonical bundle, let $n$ be an integer, and let $\cE$ be an
object of $\D(X)$ such that
\begin{itemize}
\item[--] $\cE \iso S_X\cE[-n]$;
\item[--] $\Hom_{\D(X)}(\cE, \cE) = k$;
\item[--] $\Hom_{\D(X)}(\cE, \cE[i]) = 0$ for $i<0$.
\end{itemize}
Then $\cE$ is of the form $\cO_x[i]$ for $x$ a closed point in $X$ and
$i\in \Z$.
\end{proposition} 

We would like to argue that the equivalence $\D(X)\ra \D(Y)$ induces a
bijective map between the set of (shifts of) closed points in $X$ and
the set of (shifts of) closed points in $Y$.  If we knew that $Y$ had
ample or antiample canonical class, this would follow at once from the
above proposition, but since we only know this property for $X$, we
need to work a bit harder.

Call an object satisfying the conditions of Proposition~\ref{prop:po}
a {\em point object}.  Since point objects are defined abstractly,
just in terms of categorical properties, it follows that the
equivalence $\D(X) \stackrel{\sim}{\lra} \D(Y)$ must map point objects
in $X$ to point objects in $Y$.  Since on $X$ the point objects are
all of the form $\cO_x[i]$ for some $x\in X$, $i\in \Z$, it follows
that two point objects $P,Q$ in $\D(X)$ (and hence in $\D(Y)$) must
either differ by a shift, or satisfy
\[ \Hom_{\D(X)}(P,Q[i]) = 0 \mbox{ for all } i. \]
Assume there is a point object $P$ in $Y$ which is not of the form
$\cO_y[i]$ for $y$ a closed point in $Y$ and $i\in \Z$.  Then, since
$\cO_y[i]$ is certainly a point object in any case for all $y, i$, it
follows that we must have
\[ \Hom_{\D(Y)}(P, \cO_y[i]) = 0 \mbox{ for all }y\mbox{ a closed
  point in }Y,~i\in \Z. \]
An easy argument now shows that $P$ must be zero.

We conclude that the given equivalence maps shifts of structure
sheaves of closed points in $X$ to shifts of structure sheaves of
closed points in $Y$, bijectively.

Knowing, on a variety, what the class of objects that are
isomorphic to shifts of structure sheaves of points is, enables one to
identify the class of shifts of line bundles on $X$:
\begin{proposition}[\cite{BonOrl}]
\label{prop:lbo}
Let $X$ be any smooth, quasi-projective variety.  Let $\cL$ be an
object of $\D(X)$ such that for every closed point $x$ in $X$, there
exists an integer $s$ such that
\begin{itemize}
\item[--] $\Hom_{\D(X)}(\cL, \cO_x[s]) = k$;
\item[--] $\Hom_{\D(X)}(\cL, \cO_x[i]) = 0$ if $i\neq s$.
\end{itemize}
Then $\cL$ is a shift of a line bundle on $X$.
\end{proposition}

We conclude that the given equivalence must map a line bundle on $X$
to a shift of a line bundle on $Y$.  Fix a line bundle $\cL_X$ on $X$,
and adjust the equivalence by a shift so that the image of $\cL_X$ is a
line bundle $\cL_Y$ (and not a shift thereof).

Now we can map closed points of $X$ to closed points of $Y$: these are
point objects for which the integer $s$ in Proposition~\ref{prop:lbo}
for the line bundle $\cL_X$ (or $\cL_Y$, respectively) is precisely 0.

Knowing what the points are (and not only their shifts) allows us to
identify the line bundles (and not their shifts): a line bundle is a
``line bundle object'' $\cL$ which satisfies 
\[ \Hom(\cL, \cO_x) = k \mbox{ for all } x\mbox{ a closed point in }
X. \]
This also shows that $\dim X = \dim Y$, as $S_X\cL_X[-\dim X]$ is a
line bundle, and it must be mapped by the equivalence to a line bundle
as well, but it is mapped to $S_Y\cL_Y[-\dim X]$.

Having two line bundles $\cL_1, \cL_2$, and a morphism $f$ between
them, we can now categorically identify the zero locus of this
morphism: it is the set of closed points $x\in X$ such that the
induced map
\[ \Hom(\cL_2, \cO_x) \ra \Hom(\cL_1, \cO_x) \]
is zero.  The complements of these zero loci for all line bundles and
all maps between them form a set of obviously open sets in $X$ and
hence in $Y$, and these open sets are easily seen to form a basis for
the topologies of $X$ and of $Y$, respectively.  This shows that the
map on sets of closed points between $X$ and $Y$ is continuous.

Let $\cL_X^i = S_X^i \cL_X[-i\dim X]$, and similarly for $\cL_Y$.  Since
$\omega_X$ is ample or antiample, a result of Illusie shows that the
collection of open sets obtained from morphisms between $\cL_X^i$'s
forms a basis for the topology of $X$.  Therefore the same must be
true of the $\cL_Y^i$'s, which implies that $\omega_Y$ is ample or
antiample.

Now we can let 
\[ R(X) = \bigoplus_i \Hom(\cL_X^0, \cL_X^i), \]
with the obvious induced graded ring structure.  The ring $R(X)$ is
isomorphic to the pluricanonical ring of $X$.  Since everything is
categorical it follows that defining $R(Y)$ in the same way, we have
\[ R(X) \iso R(Y). \]

Then we have isomorphisms
\[ X \iso \Proj R(X) \iso \Proj R(Y) \iso Y, \]
where the outer two isomorphisms come from the fact that $\omega_X$
and $\omega_Y$ are ample or anti-ample.

\begin{remark}
We will see in Lecture 6 that the isomorphism of pluricanonical rings
is true in general, for every equivalence of derived categories,
without any assumptions on the ampleness of the canonical class.
\end{remark}

\section{Lecture 5: Moduli problems and equivalences}

One topic that is truly central to the study of derived categories and
that we have only touched on until now is how do we define functors
between derived categories, and when are such functors equivalences?

\subsection{}
Fix smooth, projective varieties $X$ and $Y$.  It turns out that the
only way to construct exact functors $\D(X)\ra \D(Y)$ anybody knows is
that of Definition~\ref{def:intfun}: choose an object $\cE\in
\D(X\times Y)$ and consider the {\em integral transform} with {\em
kernel} $\cE$ given by
\[ \FMXY^\cE:\D(X)\ra \D(Y), \quad \FMXY^\cE(-) =
\R\pi_{Y,*}(\pi_X^*(-) \lotimes \cE). \]

The first important observation is that composing integral transforms
yields again an integral transform:
\begin{proposition}
\label{prop:compfunc}
Let $\cE\in \D(X\times Y)$, $\cF\in \D(Y\times Z)$, and define
\[ \cF\circ\cE = \R\pi_{XZ, *}(\pi_{XY}^*\cE \lotimes \pi_{YZ}^*
\cF), \]
where $\pi_{XY}$, $\pi_{YZ}$, $\pi_{XZ}$ are the projections from
$X\times Y\times Z$ to $X\times Y$, $Y\times Z$, $X\times Z$,
respectively.  Then there is a canonical isomorphism of functors
\[ \FMXZ^{\cF\circ \cE} = \FMYZ^\cF \circ \FMXY^\cE. \]
\end{proposition}

\begin{exercise}
Prove it!
\end{exercise}

\subsection{}
Another useful observation is that the functor $\FMXY^\cE$ always has
both a left and a right adjoint, respectively given by
\[ \FMYX^{\cE^\chk \otimes \omega_Y[\dim Y]}\mbox{ and }
\FMYX^{\cE^\chk \otimes \omega_X[\dim X]}, \]
where 
\[ \cE^\chk = \R\sHom_{X\times Y}(\cE, \cO_{X\times Y}).\]

\begin{exercise}
Use Serre duality to prove this.
\end{exercise}

\subsection{}
\label{subsec:exintt}
All the usual covariant functors that we know can be written as
integral transforms.  For example, if $f:X\ra Y$ is a morphism, let
$\Gamma_f\subseteq X\times Y$ be the graph of $f$.  Then
$\cO_{\Gamma_f}\in \D(X\times Y)$ will induce the $\R f_*$ functor
when viewed as a kernel $X\ra Y$, and the $\Ld f^*$ functor when
viewed as $Y\ra X$.  If $\cE\in \D(X)$, the functor $\FMXX^{\Delta_*
  \cE}$ is nothing but $-\, \lotimes \cE:\D(X) \ra \D(X)$. And so on.

\begin{exercise}
Prove these assertions.
\end{exercise}

\subsection{The elliptic curve automorphism}
We now come to the mother of all interesting examples.  It was
originally worked out by Mukai, as part of his study of derived
categories of abelian varieties.

Let $E$ be an elliptic curve with origin $P_0$, and consider the line
bundle on $E\times E$
\[ \cP = \pi_1^* \cO_E(-P_0) \otimes \pi_2^* \cO_E(-P_0) \otimes
\cO_{E\times E}(\Delta), \]
where $\Delta\subseteq E\times E$ is the diagonal.  The line bundle
$\cP$ is called the {\em normalized Poincar\'e bundle}.  Its main
property is that 
\[ \cP|_{E\times \{P\}} \iso \cP|_{\{P\}\times E} \iso \cO(P-P_0). \]

Consider the integral transform $\Phi = \Phi_{E\ra E}^\cP$ induced by
$\cP$.  

\begin{theorem}[Mukai~\cite{MukAb}]
\label{thm:mukab}
We have $\Phi \circ \Phi = \iota\circ[-1]$, where $\iota:E\ra E$ is
the negation map in the group structure on $E$.  Since
$\iota\circ[-1]$ is obviously an equivalence, it follows that $\Phi$
itself is a (non-trivial) equivalence.
\end{theorem}

To prove this, we first need a few lemmas:
\begin{lemma}[Seesaw principle]
Let $\cL, \cM$ be line bundles on the product $X\times Y$ of
projective varieties $X$ and $Y$.  Assume that $\cL_y\iso \cM_y$ for
all $y\in Y$, and that $\cL_x \iso \cM_x$ for some $x\in X$.  (Here
$\cL_y$ is defined to be $\cL|_{X\times \{y\}}$, viewed as a sheaf on
$X$, and similarly for $\cL_x$.)  Then $\cL \iso \cM$.
\end{lemma}

\begin{proof}
(Sketch.)  The condition that $\cL_y\iso \cM_y$ for all $y\in Y$
  implies that $\cL \iso \cM \otimes \pi_Y^* \cN$ for some line bundle
  $\cN$ on $Y$.  It now follows that $\cL_x \iso \cM_x \otimes \cN$,
  so $\cN$ must be trivial.
\end{proof}

For $P\in E$ a closed point, let $D_P$ be the divisor on $E\times E$
given by
\[ D_P = \{ (x, y)\in E\times E~|~ x+y = P \mbox{ in the group law of
}E\}. \]
Slightly more abstractly, if $m:E\times E \ra E$ is the group law of
$E$, then $D_P = m^{-1}(P)$.

\begin{lemma}
For $P\in E$, we have 
\[ \cO(D_P - D_{P_0}) \iso \pi_1^* \cO(P-P_0) \otimes
\pi_2^*\cO(P-P_0). \]
\end{lemma}

\begin{proof}
Apply the seesaw principle.  Fix a point $Q\in E$, and consider the
restriction of the two line bundles in question to $E\times \{Q\}$.
The right hand side obviously restricts to $\cO(P-P_0)$.  For the left
hand side, there are unique points $Q_1$ and $Q_2$ such that $Q_1+Q =
P$ and $Q_2+Q = P_0$ in the group law.  This means that $Q_1+Q-P_0\sim
P$ and $Q_2+Q-P_0\sim P_0$ (linear equivalence).  We then have
\[ \cO(D_P - D_{P_0})|_{E\times \{Q\}} \iso \cO(Q_1-Q_2) \iso \cO((Q_1
+Q -P_0) -(Q_2+Q-P_0)) \iso \cO(P-P_0). \]
Therefore the restrictions of the left and the right hand side to
every $E\times\{Q\}$ are isomorphic.  Similarly, the restrictions to
$\{Q\} \times E$ are also isomorphic for every $Q$, so by the seesaw
principle we get the desired isomorphism.
\end{proof}

\begin{lemma}
We have an isomorphism of line bundles on $E\times E\times E$
\[ \pi_{12}^* \cP \circ \pi_{23}^* \cP \iso f^* \cP, \]
where $f:E\times E\times E \ra E\times E$ is the map which takes
\[ f(x,y,z) = (x+z, y), \]
the first component being addition in the group law of $E$.
\end{lemma}

\begin{proof}
Again, we apply the seesaw principle.  It will be useful to
distinguish the copy of the elliptic curve that is mapped identically
by $f$, so we'll call it $\hat{E}$.  (The suggestion implied by the
notation is correct: everything we do is essentially the same for
abelian varieties of arbitrary dimension, but the equivalence we get
will be between $\D(A)$ and $\D(\hat{A})$, where $\hat{A}$ denotes the
dual abelian variety.)

Thus we are trying to argue that two line bundles on $E\times \hat{E}
\times E$ are isomorphic, and to do so we need to check the
isomorphism of their restrictions to fibers of the form $E\times
\{P\} \times E$, for $P\in \hat{E}$.  Over that fiber we have
\begin{align*}
\pi_{12}^* \cP |_{E\times \{P\} \times E} & \iso \pi_1^*
\cO(P-P_0) \\
\pi_{23}^* \cP |_{E\times \{P\} \times E} & \iso \pi_2^*
\cO(P-P_0) \\ 
f^* \cP |_{E\times \{P\} \times E} & \iso m^* \cO(P-P_0)  = \cO(D_P -
D_{P_0}).
\end{align*}
The previous lemma now implies that
\[ (\pi_{12}^* \cP \circ \pi_{23}^* \cP)|_{E\times \{P\} \times E} \iso f^* \cP|_{E\times \{P\} \times E} \]
for every $P\in \hat{E}$.

Now it is immediate to observe that the fibers of both line bundles
over $\{P_0\}\times \hat{E}\times \{P_0\}$ are trivial, so the seesaw
principle gives the result.
\end{proof}

\subsection{}
We can now proceed to the proof of Theorem~\ref{thm:mukab}:
\begin{proof}
By Proposition~\ref{prop:compfunc}, the composition $\Phi\circ \Phi$
is given by $\Phi_{E\ra E}^{\cP\circ \cP}$.  (Strictly speaking, the
second $\cP$ ought to be $\tau^* \cP$, where $\tau:E\times E\ra
E\times E$ is the involution interchanging the two factors; but it is
obvious that $\tau^*\cP \iso \cP$ from the definition of $\cP$.)  We
have, by the previous Lemma,
\[ \cP\circ \cP = \R \pi_{13,*} (\pi_{12}^* \cP \otimes \pi_{23}^*
\cP) \iso \R \pi_{13,*} f^* \cP. \]
From the flat-base change formula for the diagram
\[ 
\begin{diagram}[height=2em,width=2em,labelstyle=\scriptstyle]
E\times \hat{E} \times E & \rTo^f & E\times \hat{E} \\
\dTo^{\pi_{13}} & & \dTo^{\pi_1} \\
E\times E & \rTo^m & E
\end{diagram}
\]
we conclude that
\[ \R \pi_{13,*} f^* \cP \iso m^* \R\pi_{1,*} \cP. \]

The crucial computation (exercise below) is that $\pi_{1,*} \cP = 0$
and $\R^1 \pi_{1,*} \cP = \cO_{P_0}$, the structure sheaf of the
origin.  Knowing this, it follows that
\[ \R\pi_{1,*} \cP = \cO_{P_0}[-1], \]
and therefore
\[ m^* \R\pi_{1,*} \cP = \cO_{\Gamma_\iota}[-1]. \] 
From the basic examples in~(\ref{subsec:exintt}) we get the result of
the theorem.
\end{proof}

\begin{exercise}
Prove that $\pi_{1,*} \cP = 0$ and $\R^1 \pi_{1,*} \cP = \cO_{P_0}$
where $\pi_1:E\times \hat{E} \ra E$ is the projection, and $\cP$ is
the Poincar\'e bundle.
\end{exercise}

\noindent \textbf{Hint.} Prove that
\[ h^i(\hat{E}, \cP_e)  = \left \{ \begin{array}{ll} 0 & \mbox{ for }
  e\neq P_0 \\ 1 & \mbox{ for } e=P_0 \end{array}\right . \\
\]
for $i=0,1$ and any $e\in E$.  Using~\cite[Ex. III.11.8]{HarAG} it
follows that $\R^i \pi_{1,*} \cP$ must be supported at $P_0$ for
$i=0,1$.  Now argue that $H^0(E\times \hat{E}, \cP) =0$, so $\R^0
\pi_{1,*} \cP = 0$.  Use Riemann-Roch on the surface $E\times \hat{E}$
to argue that $\chi(\cP) = -1$, so $h^1(E\times \hat{E}, \cP) = 1$,
and finish using the Leray spectral sequence.

\subsection{}
Observe that in this example we have
\[ \Phi(\omega_E) = \Phi(\cO_E) = \cO_{P_0} \neq \omega_E, \]
so we have an equivalence which does not map $\omega_X$ to $\omega_Y$.

\subsection{}
We worked out Mukai's example in full detail to have an explicit
calculation around, but since 1981 (the year of Mukai's result) people
have developed much more powerful techniques for checking when an
integral transform is an equivalence.  In fact, integral transforms
that are equivalences have a special name, {\em Fourier-Mukai
  transforms}.  The reason for the name Fourier is the following
imprecise analogy with the Fourier transform.

Think of objects of $\D(X)$ as being the analogue of smooth functions
on the circle $S^1$.  The analogue of the usual metric on the space of
functions is given by the (collection of) pairings
\[ \langle \cE, \cF\rangle_i = \dim \Hom_{\D(X)}(\cE, \cF[i]). \]
Under this analogy, an equivalence $\Phi$ is analogous to an isometry
between two spaces of functions, which is precisely what the usual
Fourier transform is.  In fact, the analogy is even more precise: the
usual Fourier transform is given by integrating (=push-forward along
second component) the product (=tensor product) of the original
function (=pull-back of the object of $\D(X)$) with a kernel (=the
kernel object on $\D(X\times Y)$).

The trick to proving that the usual Fourier transform is an isometry
is to find an orthonormal basis (in an appropriate sense) for the
space of functions, and to prove that the image of this orthonormal
basis is another orthonormal basis.  It turns out that the same exact
trick can be used to check when an integral transform is an
equivalence: an ``orthonormal basis'' for the collection of objects in
$\D(X)$ is given by $\{\cO_x\}_{x\in X}$, in the sense that
\[ \Ext^i(\cO_x, \cO_y) = \left \{ \begin{array}{ll}
0 & \mbox{ if } x\neq y\mbox{ or } i\not\in [0,\dim X] \\
\C & \mbox{ if } x = y \mbox{ and } i = 0.
\end{array} \right . \]
It turns out that to check that an integral transform $\Phi$ is fully
faithful (=preserves inner products), it will suffice to check that the
image of the orthonormal basis $\{\cO_x\}$ under $\Phi$ is again
orthonormal, in the sense of the above equations.
\begin{theorem}[Mukai, Bondal-Orlov, Bridgeland]
\label{thm:equivcrit}
Let $X, Y$ be smooth, projective varieties, and let $\Phi:\D(X)\ra
\D(Y)$ be an integral transform.  For $x\in X$ let $\cP_x =
\Phi(\cO_x)$.  Then $\Phi$ is fully faithful if, and only if, we have
\[ \Hom_{\D(Y)}(\cP_x, \cP_y[i]) = \left \{ \begin{array}{ll}
0 & \mbox{ if } x\neq y\mbox{ or } i\not\in [0,\dim X] \\
\C & \mbox{ if } x = y \mbox{ and } i = 0.
\end{array} \right . \]

Furthermore, $\Phi$ is an equivalence if, and only if, 
\[ \cP_x\otimes \omega_Y \iso \cP_x \]
for all $x\in X$.
\end{theorem}

(It is easy to see that the condition $\cP_x\otimes \omega_X \iso
\cP_x$ must be satisfied: this is because $\cO_x \otimes \omega_X \iso
\cO_x$, and equivalences commute with Serre functors.)

\begin{exercise}
Check that Mukai's original equivalence is indeed an equivalence using
the above criterion.
\end{exercise}

\subsection{Moduli point of view}
The above results lend a new perspective to what an integral transform
is: it is a way to view the objects $\{\cP_x\}_{x\in X}$ as a
collection of (complexes of) sheaves on $Y$, parametrized by the
points of $X$.  This immediately brings to mind {\em moduli spaces of
sheaves}: a space $M$ is called a moduli space of (say, stable)
sheaves on $X$ if the points of $M$ are in 1-1 correspondence with the
set of isomorphism classes of sheaves with some very specific
properties (e.g., stable, with fixed Chern classes, etc.)

The most concrete construction of this sort can be found again in
Mukai's example: the points of the dual elliptic curve $\hat{E}$
(which happens to be isomorphic to $E$) parametrize all the line
bundles of degree 0 on $E$ (via the correspondence $P\in \hat{E} = E
\mapsto \cO_E(P-P_0)$).  

It turns out that stable sheaves are also good in general (to some
extent) from the orthogonality point of view: a standard theorem
about stable sheaves implies that if $\cE, \cE'$ are stable sheaves
(with respect to some polarization) on a projective variety $X$, and
the Hilbert polynomials of $\cE$, $\cE'$ are the same, then 
\[ \Hom_X(\cE, \cE') = \left \{ \begin{array}{ll} 0 & 
\mbox{ if } \cE \neq \cE' \\ k & \mbox{ if } \cE =
\cE'.\end{array}\right . \]

\subsection{}
As an example of an application of this, consider the following result
of Mukai:
\begin{theorem}
Let $X$ be an abelian or K3 surface, and let $M$ be a moduli space of
stable sheaves such that
\begin{enumerate}
\item $M$ is compact and 2-dimensional;
\item the Chern classes of the sheaves being parametrized are such
  that their relative Euler characteristic is equal to 2;
\item the moduli problem is fine, i.e., there exists a universal sheaf
  $\cU$ (the analogue of the Poincar\'e bundle from before)
\end{enumerate}
Then $\FMMX^\cU:\D(M) \ra \D(X)$ is an equivalence.
\end{theorem}

(Mukai's result is in fact much more powerful: it gives complete
numerical criteria for conditions 1), 2) and 3) to be satisfied.  And
Mukai is able to give a complete description of $M$, in terms of the
initial K3 surface and numerical data used to specify the moduli
problem, via the Torelli theorem.)

\begin{proof}
We only need to check the conditions of Theorem~\ref{thm:equivcrit}.
The objects $\cP_x$ of the theorem become $\cU_m$, for $m\in M$, i.e.,
precisely the stable sheaves parametrized by $M$.  Since $X$ and $M$
have dimension 2 and $\cU_m$ are sheaves (as opposed to complexes),
the condition that $\Ext^i_X(\cU_m, \cU_m') = 0$ if $i\not\in [0,2]$
is automatic.  The fact that $\Hom_X(\cU_m, \cU_m) = k$ follows from
the fact that $\cU_m$ is stable, hence simple for all $m$.  If $m\neq
m'$, we have $\Hom_X(\cU_m, \cU_{m'}) = 0$ from stability; Serre
duality on $X$ implies that $\Ext^2_X(\cU_m, \cU_{m'}) = 0$; and the
fact that the relative Euler characteristic of $\cU_m, \cU_{m'}$ is
zero implies that $\Ext_X^1(\cU_m, \cU_{m'}) = 0$.
\end{proof}

\section{Lecture 6: Hochschild homology and cohomology}

\subsection{}
The derived category is a rather involved object, so one of the things
one would like to have is some more manageable invariants.  More
explicitly, it would be highly desirable to have some numerical or
vector space invariants (analogues of Betti numbers or
homology/cohomology) that one could associate to a space, and that
would be invariant under equivalences of derived categories.

One reason this is interesting is because in the past few years there
have been many powerful results that have been phrased as equivalences
between seemingly unrelated derived categories.  Perhaps one of the
more striking is the Bridgeland-King-Reid equivalence:
\begin{theorem}[Bridgeland-King-Reid~\cite{BKR}]
Let $M$ be a quasi-projective variety endowed with the action of a
finite group $G$ such that the canonical class of $M$ is locally
trivial as a $G$-sheaf.  Let $X = M/G$, and let $Y=\GHilb(M)$.  If
$\dim Y\times_X Y = \dim M + 1$, then $Y$ is a crepant resolution of
$X$ and there is a natural equivalence of derived categories
\[ \D(Y) \iso \D([M/G]), \]
where the latter is the derived category of $G$-equivariant sheaves on
$M$.
\end{theorem}

One would like to be able to extract numerical information out of
this; for example, one would like to be able to compute the
topological Euler characteristic of the crepant resolution $Y$ when
$M$ and the action of $G$ are well understood.

\subsection{}
The invariants we will be talking about today all arise from the
following observation: there are three natural, intrinsic functors
$\D(X)\ra \D(X)$ for every smooth, projective $X$ -- namely, the
identity functor, the Serre functor, and the shift functor.  Any
equivalence $\D(X)\ra \D(Y)$ must commute with these, {\em and must
  take natural transformations between them to corresponding natural
  transformations}.  It is this second part of the observation that
leads one to the following definition:
\begin{definition}
Let $X$ be a smooth, projective variety.  Define the Hochschild
cohomology of $X$ to be given by
\[ HH^*(X) = \Ext_{X\times X}^*(\cO_\Delta, \cO_\Delta), \]
and the Hochschild homology of $X$ by
\[ HH_*(X) = \Ext_{X\times X}^{-*}(\Delta_! \cO_X, \cO_\Delta). \]
Here, $\Delta_!:\D(X) \ra \D(X\times X)$ is the left adjoint of $\Ld
\Delta^*$, defined using the techniques of
Theorem~\ref{thm:leftrightadj}.
\end{definition}

\subsection{}
Before we turn to proving that these vector spaces are invariants of
the derived category of $X$, let us compute them when $X$ is a curve
$C$.  To compute $\Ext^i(\cO_\Delta, \cO_\Delta)$ we'll use the
local-to-global spectral sequence, so we need to compute first
$\sExt^i(\cO_\Delta, \cO_\Delta)$.  Since $C$ is a curve, $\Delta$ is
a divisor in $C\times C$, so we have the locally free resolution
\[ 0 \ra \cO_{C\times C}(-\Delta) \ra \cO_{C\times C} \ra \cO_\Delta
\ra 0. \]
Taking $\sHom(-,\cO_\Delta)$ yields the following complex that
computes $\sExt^*(\cO_\Delta, \cO_\Delta)$:
\[ 0 \ra \cO_\Delta \stackrel{0}{\lra} \cO(-\Delta)|_\Delta =
\cN_{\Delta/C\times C} \ra 0, \] 
where the second term is the normal bundle of $\Delta$ in $C\times C$
(viewed as a sheaf on $C\times C$ by push-forward), i.e., $\Delta_*
T_C$.  The spectral sequence now looks like
\[
\begin{diagram}[height=2em,width=2em]
\vdots & &\vdots &&\vdots &&\\
\uTo & & \uTo && \uTo && \\
0 & \rTo & 0 & \rTo & 0 & \rTo & \cdots\\
\uTo & & \uTo && \uTo && \\
H^0(C, T_C) & \rTo & H^1(C, T_C) & \rTo & 0 & \rTo & \cdots\\
\uTo & & \uTo && \uTo &&  \\
H^0(C, \cO_C) & \rTo & H^1(C, \cO_C) & \rTo & 0 & \rTo & \cdots
\end{diagram}
\]
and all the maps are zero.  We conclude that we have
\begin{align*}
HH^0(C) & = H^0(C, \cO_C) \\
HH^1(C) & = H^1(C, \cO_C) \oplus H^0(C, T_C) \\
HH^2(C) & = H^1(C, T_C).
\end{align*}

\subsection{}
\label{subsec:hhcurve}
Let us now move to homology: Theorem~\ref{thm:leftrightadj} shows that
$\Delta_! \cO_X \iso \Delta_* \omega_X^{-1}[-\dim X]$.  Running through
the same calculations as before, we find that $HH_i$ will be non-zero
for $-1\leq i\leq 1$ and we have
\begin{align*}
HH_{-1}(C) & = H^1(C, \cO_C) \\
HH^0(C) & = H^0(C, \cO_C) \oplus H^1(C, \Omega_C) \\
HH^1(C) & = H^0(C,\Omega_C).
\end{align*}

Based on these calculations, it should come as no surprise that we
have
\begin{theorem}[Hochschild-Kostant-Rosenberg~\cite{HKR},
    Swan~\cite{Swa}, Kontsevich~\cite{KonDQ}, Yekutieli~\cite{Yek}]
\label{thm:hkr}
If $X$ is a smooth, quasi-projective variety over $\C$, then
\begin{align*}
HH^i(X) & \iso \bigoplus_{p+q=i} H^p(X, \wedge^q T_X) \\
HH_i(X) & \iso \bigoplus_{q-p=i} H^p(X, \Omega_X^q).
\end{align*}
\end{theorem}

Especially in the case of homology, this theorem gives the connection
with the more topological Hodge and Betti numbers: $HH_i(X)$ is the
sum of the vector spaces on the $i$-th {\em vertical} of the Hodge
diamond of $X$, counting from the center.

\subsection{}
What does the definition of $HH^*(X)$, $HH_*(X)$ have to do with the
three intrinsic functors we originally listed and the natural
transformations between them?  We have
\[ HH^i(X) = \Ext_{X\times X}^i(\cO_\Delta, \cO_\Delta) =
\Hom_{\D(X\times X)}(\cO_\Delta, \cO_\Delta[i]). \]
Now observe that, as a kernel, $\cO_\Delta$ corresponds to the
identity functor $\D(X)\ra \D(X)$, and $\cO_\Delta[i]$ corresponds to
the ``shift by $i$'' functor.  It is easy to see that
given a morphism (in $\D(X\times X)$) $\cO_\Delta \ra \cO_\Delta[i]$,
one gets for every $\cF\in\D(X)$ a map 
\[ \cF  = \FMXX^{\cO_\Delta}(\cF) \ra \FMXX^{\cO_\Delta[i]}(\cF) =
\cF[i], \] 
and, furthermore, these maps form a natural transformation $\Id
\Rightarrow [i]$.

\subsection{}
At this point it is easy to state a potential theorem:
\begin{conjecture}
\label{conj:prod}
The functor
\[ \Phi: \D(X\times Y) \ra \ExFun(\D(X), \D(Y)) \]
is an equivalence, where the objects of the right hand side are exact
functors $\D(X)\ra \D(Y)$, morphisms between them are natural
transformations that commute with shifts, and $\Phi$ takes an object
$\cE$ of $\D(X\times Y)$ to the integral transformation $\FMXY^\cE$,
and morphisms between objects to the obvious natural transformations.
\end{conjecture}
If this conjecture were true, then it would be obvious that our
definition of Hochschild cohomology is completely intrinsic to
$\D(X)$: $HH^i(X)$ would equal the space of natural transformations
that commute with shifts between the identity functor and the
shift-by-$i$ functor.

However, tough luck: the conjecture is most definitely false, as the
following example demonstrates.
\begin{example}
If $E$ is an elliptic curve we have $HH^2(E) = k$, but any natural
transformation between the identity functor and the shift-by-2 functor
is zero. 

The calculation that $HH^2(E) = k$ is~(\ref{subsec:hhcurve}).  The
fact that any natural transformation $\Id\Rightarrow [2]$ is zero
follows from the fact that on any curve $\Ext^2_E(\cF, \cF) = 0$ for
any sheaf $\cF$, and the fact that any object in the derived category
of a curve is isomorphic to a direct sum of its cohomology sheaves,
which in turn is a consequence of the fact that the cohomological
dimension of the curve is 1. (\textbf{Exercise:} Fill in the details.)
\end{example}

\subsection{}
The failure of Conjecture~\ref{conj:prod} is the main reason it is
believed the derived category should be replaced by a more powerful
technical tool; possible choices here are differential-graded (dg)
categories, or $\cA_\infty$-categories (these two classes are
equivalent to each other).  For dg-categories, the corresponding
statement to Conjecture~\ref{conj:prod} has recently been proved by
Toen~\cite{Toen}. 

Because of this failure, we cannot simply define $HH^*(X)$ as
``natural transformations...''  The definition that we gave is what
the dg approach would define.  But now we have to work harder to prove
that cohomology (and homology) are derived category invariants.

To prove this, we first need a powerful result of Orlov~\cite{Orl}:
\begin{theorem}
\label{thm:equivkern}
Let $\Phi:\D(X)\ra \D(Y)$ be an equivalence between the derived
categories of smooth, projective varieties $X$ and $Y$.  Then there
exists a kernel $\cE\in\D(X\times Y)$ such that 
\[ \Phi \iso \FMXY^\cE, \]
and any two such kernels are connected by a canonical isomorphism.
\end{theorem}

Note what this theorem says: even though we do not know that the image
of the functor $\Phi$ from Conjecture~\ref{conj:prod} is essentially
surjective, Orlov's result says that it is essentially surjective if
we restrict our attention to equivalences (and kernels giving them).

\begin{proposition}
\label{prop:conj}
Let $X$, $Y$, $Z$, $W$ be smooth projective varieties, and let $\cE\in
\D(X\times Y)$, $\cG\in \D(Z\times W)$ be kernels of integral
transforms $\FMXY^\cE$, $\Phi_{Z\ra W}^\cG$.  Consider 
\[ \cH = \cE \boxtimes \cG \in \D(X\times Y\times Z\times W), \]
and view it as a kernel for an integral transform 
\[ C = \Phi_{Y\times Z\ra X\times W}^\cH: \D(Y\times Z) \ra
\D(X\times W). \]
Then we have, for $\cF\in \D(Y\times Z)$,
\[ \Phi_{X\ra W}^{C(\cF)} \iso \Phi_{Z\ra W}^\cG \circ \FMYZ^\cF \circ
\FMXY^\cE. \]
\end{proposition}

\begin{exercise}
Use the appropriate base change theorem to prove this.
\end{exercise}

We also need the following theorem of Orlov:
\begin{theorem}[\cite{OrlAb}]
\label{thm:equivprod}
Assume $X, X', Y, Y'$ are smooth, projective varieties, and $\cE\in
\D(X\times Y)$, $\cE'\in\D(X'\times Y')$ are kernels that induce
equivalences $\FMXY^\cE$, $\Phi_{X'\ra Y'}^{\cE'}$.  Then 
\[ \cE\boxtimes\cE' = \pi_{13}^* \cE \lotimes \pi_{24}^* \cE'\in
\D(X\times X'\times Y\times Y') \]
induces an equivalence
\[ \Phi_{X\times X'\ra Y\times Y'}^{\cE\boxtimes \cE'}:\D(X\times X')
\ra \D(Y\times Y'). \]
\end{theorem}

\begin{exercise}
Use the criterion for an equivalence to prove this directly.  Or
alternatively, use the ideas of Proposition~\ref{prop:conj} to
construct an inverse directly.
\end{exercise}

\begin{proposition}
\label{prop:equivdual}
Assume $\cE\in \D(X\times Y)$ induces an equivalence $\D(X)\ra
\D(Y)$.  Then $\FMXY^{\cE^\chk}:\D(X) \ra \D(Y)$ is also an
equivalence, where
\[ \cE^\chk = \R \sHom_{X\times Y}(\cE, \cO_{X\times Y}). \]
\end{proposition}

\begin{proof}
Use the fact that on a smooth, projective variety we have
\[ \Hom(\cE, \cF) \iso \Hom(\cF^\chk, \cE^\chk) \]
and check the orthogonality conditions of Theorem~\ref{thm:equivcrit}.
\end{proof}

\subsection{}
We can now turn to the proof of the invariance of homology and
cohomology under derived equivalences:
\begin{theorem}[\cite{CalHH1}]
Let $X$ and $Y$ be smooth, projective varieties, and for $i\in \Z$ let
\[ \omega_{\Delta_X}^{\otimes i} = \Delta_* \omega_X^{\otimes i}, \]
and similarly for $Y$.  Let $\Phi:\D(X) \ra \D(Y)$ be an equivalence.
Then there exists an induced equivalence
\[ \bar{\Phi}:\D(X\times X) \ra \D(Y\times Y) \]
which satisfies
\[ \bar{\Phi}(\omega_{\Delta_X}^{\otimes i}[j]) \iso
\omega_{\Delta_Y}^{\otimes i}[j]. \]

Therefore we have induced isomorphisms
\[ HH^*(X) \iso HH^*(Y),\quad HH_*(X) \iso HH_*(Y),\quad R(X) \iso
R(Y), \]
where $R(X)$ denotes the pluricanonical ring of $X$.
\end{theorem}

\begin{proof}
(Parts of this theorem also proven independently
in~\cite{GolLunOrl},~\cite{OrlAb}.)  By Theorem~\ref{thm:equivkern}
$\Phi$ is given by an integral transform; let $\cE\in\D(X\times Y)$ be
the kernel.  Let $\cE^* = \cE^\chk \otimes \pi_Y^* \omega_Y[\dim Y]$,
so that $\Phi^\dagger = \FMYX^{\cE^*}$ is an adjoint (and, hence, an
inverse) to $\FMXY^\cE$.  (Since we are dealing with an equivalence,
the left and right adjoints are naturally isomorphic to each other and
they give inverses to $\Phi$.)  By Proposition~\ref{prop:equivdual},
$\FMXY^{\cE^*}$ is also an equivalence, so by
Theorem~\ref{thm:equivprod} we get an equivalence
\[ \bar{\Phi} = \Phi_{X\times X\ra Y\times Y}^{\cE\boxtimes \cE^*}:
\D(X\times X) \ra \D(Y\times Y). \]

We claim that $\bar{\Phi}$ is the desired equivalence.  Indeed, by
Proposition~\ref{prop:conj} we know that if $\cF$ is any kernel on
$X\times X$, then $\bar{\Phi}(\cF)$, as a kernel on $Y\times Y$, will
give a functor isomorphic to 
\[ \Phi \circ \FMXX^\cF \circ \Phi^\dagger : \D(Y) \ra \D(Y). \]
If $\cF = \omega_{\Delta_X}^{\otimes i}[j]$ then $\FMXX^\cF$ is the
functor of tensoring with $\omega_X^{\otimes i}$ and shifting by $j$.
Let $\cG = \omega_{\Delta_Y}^{\otimes i}[j]$ induce the corresponding
functor on $Y$.  Since $\FMXX^\cF$ is just a shift of a power of the
Serre functor, it commutes with an equivalence:
\[ \Phi \circ \FMXX^\cF \circ \Phi^\dagger \iso \FMYY^\cG \circ
\Phi\circ \Phi^\dagger \iso \FMYY^\cG. \]
We conclude that $\bar{\Phi}(\cF)$ gives, as a kernel, an equivalence
$\D(Y)\ra \D(Y)$ isomorphic to $\FMYY^\cG$.  By Orlov's
theorem~\ref{thm:equivkern}, we have an induced isomorphism
\[ \bar{\Phi}(\omega_{\Delta_X}^{\otimes i}[j]) \iso
\omega_{\Delta_Y}^{\otimes i}[j]. \]

This immediately implies the claims about the invariance of $HH^*$ and
$HH_*$.  The claim about the invariance of the pluricanonical rings
follows from the easy observation that 
\[ \Hom_{\D(X\times X)}(\Delta_* \cF, \Delta_* \cG) \iso
\Hom_{\D(X)}(\Ld \Delta ^* \Delta_* \cF, \cG) \iso \Hom_X(\cF, \cG) \]
when $\cF, \cG$ are {\em sheaves} on a space $X$ (even though $\Ld
\Delta^*\Delta_* \cF$ has more non-zero terms than $\cF$, morphisms in
the derived category can never go to the right).

It is worth observing at this point that $HH^*(X)$ is a graded ring,
and that the isomorphism $HH^*(X) \iso HH^*(Y)$ we have constructed
respects the graded ring structures.
\end{proof}

\subsection{}
We conclude by listing, without proof, several other fundamental
results about Hochschild homology and cohomology:
\begin{theorem}[\cite{CalHH1}]
Hochschild homology is functorial with respect to integral
transforms.  Namely, to every integral transform $\Phi:\D(X)\ra \D(Y)$
there is an associated map of graded vector spaces $\Phi_*:HH_*(X) \ra
HH_*(Y)$ such that $(\Phi\circ \Psi)_* = \Phi_* \circ \Phi_*$ and
$(\Id)_* = \Id$.
\end{theorem}

\begin{theorem}[\cite{CalHH1}]
\label{thm:adj}
There is a natural, non-degenerate, graded bilinear pairing 
\[ HH_*(X) \otimes HH_*(X) \ra \C, \]
called the generalized Mukai pairing.  If $\Phi \adjoint \Psi$, then
$\Phi_*$ is left adjoint to $\Psi_*$ with respect to the pairings on
$HH_*(X)$ and $HH_*(Y)$, respectively.
\end{theorem}

\subsection{}
An (almost) immediate consequence of the above two properties is a
form of the Hirzebruch-Riemann-Roch theorem.  To talk about it, note
first that there is a Chern class map
\[ \ch:\D(X) \ra HH_0(X), \]
constructed as follows: we have $HH_0(\pt) \iso \C$, and there is a
distinguished element $\bone\in HH_0(\pt)$ corresponding to the
identity.  An element $\cE\in\D(X) = \D(\pt\times X)$ induces an
integral transform 
\[ \Phi_{\pt\ra X}^\cE:\D(\pt) \ra \D(X), \]
and as such there is an induced map
\[ (\Phi_{\pt\ra X}^\cE)_*:HH_0(\pt) \ra HH_0(X). \]
If we define 
\[ \ch(\cE) = (\Phi_{\pt\ra X}^\cE)_* \bone \in HH_0(X) \]
it can be proven~(\cite{CalHH2}) that the image of $\ch(\cE)$ under the isomorphism of
Theorem~\ref{thm:hkr} yields the usual Chern character 
\[ \ch(\cE) \in \bigoplus_p H^p(X, \Omega^p_X). \]

An easy application of Theorem~\ref{thm:adj} yields the formal
Hirzebruch-Riemann-Roch theorem
\begin{theorem}[\cite{CalHH1}]
For $\cE, \cF\in \D(X)$ we have
\[ \langle \ch(\cE), \ch(\cF) \rangle = \chi(\cE, \cF) = \sum_i (-1)^i
\dim \Ext^i_X(\cE, \cF). \]
\end{theorem}

\subsection{}
To return to applications of the Hochschild theory to orbifolds, we have
\begin{theorem}[\cite{CalHH1},\cite{Bar}]
Hochschild homology can equally well be defined for an orbifold
(Deligne-Mumford stack), and for a global quotient we have
\[ HH_*([X/G]) \iso \left ( \bigoplus_{g\in G} HH_*(X^g) \right
)^G. \]
Here, $X^g$ is the locus of $X$ fixed by the element $g\in G$, and the
superscript $G$ denotes the fixed part of the corresponding vector
space.
\end{theorem}

This theorem can be applied, in conjunction with the derived
equivalence invariance of Hochschild homology and with the
Bridgeland-King-Reid result to argue that, in the cases where the BKR
theorem can be applied, we have~(\cite{DHVW})
\[ \chi(Y) = \sum_{[g]\in\Cl(G)} \chi(X^g/N_g^G), \]
where the sum is over conjugacy classes of elements in $G$, $Y$ is a
crepant resolution of the singularities of $X/G$, and $N_g^G$ denotes
the normalizer of $g$ in $G$.  Thus, we arrive at what we aimed for in
the beginning, namely to extract out of derived equivalences numerical
information about the underlying spaces.

\end{document}